\documentclass{article} 

\usepackage{amsmath,amsfonts,amsthm,amssymb}
\usepackage{tikz}
\usepackage[hyphens]{url}
\usepackage{hyperref}
\usepackage{pgfplots}
\usepackage{pgflibraryarrows, booktabs, mathtools}
\usepackage{pgflibrarysnakes}
\usepackage{pgfpages}
\usetikzlibrary{patterns}
\usepackage{pgfplotstable}
\usepackage[singlelinecheck=off, justification=justified]{caption}
\usepackage{subcaption}
\usepackage{xspace}
\usepackage[margin =1in]{geometry}

\newcommand\R{{\mathbb R}}

\newcommand{\dmax}{\displaystyle \max}
\newcommand{\qBB}{{\tt quadprogBB}\xspace}
\newcommand{\qIP}{{\tt quadprogIP}\xspace}
\renewcommand{\S}{\mathcal{S}}
\newcommand{\KKTQCQP}{$\operatorname{QP_{\rm KKT}}$}
\newcommand{\hoff}{\mathcal{H}_{A,b}}
\newcommand{\hoffgen}{\mathcal{H}_{A,b}}

\newcommand{\hoffpgen}{\overline{\mathcal{H}}_{A}}
\newcommand{\BARON}{{\tt BARON}\xspace}
\newcommand{\cplex}{{\tt CPLEX}\xspace}
\newcommand{\tr}{^\intercal}

\DeclareMathOperator{\st}{s.t.}

\mathtoolsset{showonlyrefs}

\usepackage{natbib}
 \bibpunct[, ]{(}{)}{,}{a}{}{,}%
 \def\newblock{\ }%
 \def\myplotscale{0.75}

 \usepackage{authblk}

 \newtheorem{theorem}{Theorem}

\newtheorem{example}{Example}
\newtheorem{definition}{Definition}
\newtheorem{proposition}{Proposition}
\newtheorem{remark}{Remark}

\begin{document}




\title{Globally solving Non-Convex Quadratic Programs via Linear Integer Programming techniques}

\author[1]{Wei Xia\thanks{wex213@lehigh.edu}}
\author[2]{Juan Vera\thanks{j.c.veralizcano@uvt.nl}}
\author[1]{Luis F. Zuluaga\thanks{luis.zuluaga@lehigh.edu, Corresponding author.}}
\affil[1]{Department of Industrial and Systems Engineering, Lehigh University \authorcr H.S. Mohler Laboratory, 200 West Packer Avenue, Bethlehem, PA 18015}
\affil[2]{Department of Econometrics and Operations Research, Tilburg University \authorcr 5000 LE, Tilburg, The Netherlands}

\renewcommand\Authands{ and }

\maketitle

\abstract{%
Quadratic programming (QP) is a well-studied fundamental NP-hard optimization problem
which optimizes a quadratic objective over a set of linear constraints.
In this paper, we reformulate QPs as a mixed-integer linear problem (MILP).
This is done via the reformulation of QP as a linear complementary problem,
and the use of binary variables and big-M constraints, to model the complementary constraints. To obtain such reformulation, we show how to impose bounds on the dual variables without eliminating
all the (globally) optimal primal solutions; using some fundamental results on the solution of perturbed linear systems.

Reformulating non-convex QPs as MILPs provides an advantageous way to obtain global solutions as it
allows the use of current state-of-the-art MILP solvers. To illustrate this, we compare the performance of our solution approach, labeled \qIP, with the
current benchmark global QP solver  \qBB, as well as with \BARON, one of the leading non-linear programming (NLP) solvers, and \cplex's non-convex QP solver, on a large variety of QP test instances.
In practice, \qIP is shown to typically outperform by orders of magnitude \qBB, \BARON, and \cplex on standard QPs. For general QPs, \qIP outperforms \qBB, outperforms \BARON in most instances, while \cplex performs the best on these instances. For box-constrained QPs, \qIP has a comparable performance to \qBB and \BARON in small- to medium-scale instances, but is outperformed by these solvers on large-scale instances; while \cplex performs the best on box-constrained QP instances. Also, unlike \qBB, the solution approach proposed here is able to solve QP instances whose dual feasible set is unbounded.
 The {\tt MATLAB} code, called \qIP, and the instances used to perform these numerical experiments are publicly available at \url{https://github.com/xiawei918/quadprogIP}.
}%


%


\section{Introduction}
\label{intro}
{\em Quadratic programmming} (QP), is a fundamental optimization problem with a quadratic objective and
linear constraints.  QP  is NP-hard~\citep[see, e.g.,][and the references therein]{PardV91},
however, when the objective is convex, QP can be globally solved (within a predetermined precision $\epsilon > 0$) in polynomial time
via {\em interior-point methods}~\citep[see, e.g.,][]{Rene01}.
Here, the focus is on obtaining global solutions for non-convex QP.
QP is arguably the most basic instance of a (non-convex), {\em non-linear program} (NLP).
At a fundamental level, the complexity of globally solving QP lies in the fact that multiple of its local optimal solutions may not necessarily be global optimal solutions~\citep[see, e.g.][]{Bert99}.

QPs commonly arise
in applications in engineering, pure and social sciences, finance, and economics~\citep[see, e.g.,][]{HorsPT00}.
As a result,
there has been extensive work on studying how to obtain global solutions of QPs using both NLP techniques \citep[see,][for surveys in this area]{gould2003, gao2004},
and convex optimization techniques \citep[consider, e.g.,][among many others]{nesterov1998, kim2001, kim2003, chen2012}.

In this paper, we reformulate QPs as a mixed-integer linear problem (MILP). This provides an advantageous way to obtain global solutions as it
allows the use of current state-of-the-art MILP solvers.
Moreover, the numerical experiments of Section~\ref{sec:3} show that
a basic implementation of the proposed algorithm, which we refer to as {\tt quadprogIP}, typically outperform by orders of magnitude \qBB, \BARON, and \cplex on standard QPs. For general QPs, \qIP outperforms \qBB, outperforms \BARON in most instances, while \cplex performs the best on these instances. For box-constrained QPs, \qIP has a comparable performance to \qBB and \BARON in small- to medium-scale instances, but is outperformed by these solvers on large-scale instances; while \cplex performs the best on box-constrained QP instances. 

Unlike \qBB, the solution approach proposed here is able to solve QP instances whose dual feasible set is unbounded.
The {\tt MATLAB} code and the instances used
to perform these numerical experiments are available at \url{https://github.com/xiawei918/quadprogIP}.

To obtain the proposed MILP-reformulation (see Sec.~\ref{sec:SolApproach}), the QP's KKT conditions are used
to reformulate the QP as a {\em linear complementarity problem} (LCP). 
In this reformulation, the complexity of the problem is captured by the complementarity constraints.
The {\em KKT-branching approach} \citep{burer2009}, which consist on branching on this complementarity constraints, is not useful on this reformulation of the problem,
as the underlying linear relaxations at the root node of the KKT-branching tree are (under mild assumptions) unbounded
 \citep[][Cor. 2.3]{burer2009}.
  Another alternative, namely, reformulating the complementarity constraints using binary variables and big M constraints,
 requires the knowledge of bounds on the problem's KKT multipliers, which in general are unbounded \citep[cf.,][Sec. 6.1 and 6.2]{Pang12}.
To directly use MILP solvers for the solution of the QP, we overcome this requirement by restricting our attention to a subset of optimal KKT points.
We show (Theorem \ref{thm:boundDual}) that it is possible to impose bounds on the dual variables without eliminating all the (globally) optimal primal solutions.
Our results are based on  fundamental results on the approximate solution of systems of linear equations \cite[e.g.,][]{GuleHR95,  Mang81}.
 One advantage of the proposed methodology is that unlike previous related work,
 the convergence of the MILP-based approach to
the QP's global optimal solution in finite time
 follows in straightforward fashion (see Sec.~\ref{sec:graphs}). Also, the methodology can be applied to QPs without the need for assumptions
 on the relative interior of its feasible set (see Sec.~\ref{sec:comphauff} for details).

Before stating the results described above, we end this section with a short review of both NLP and convex optimization techniques for
the global solution of QP's.
Using NLP techniques,
\cite{van1999} proposed an interior-point algorithm for NLPs (thus, applicable for QPs), which
is an extension of the interior-point methods for linear and convex optimization problems \citep[cf.,][]{Rene01}. \cite{flou1990} proposed an algorithm which globally solves certain
classes of NLPs  by decomposing the problem based on an appropriate partition of its decision variables.
The work of  \cite{BeloLLMW09} and \cite{TamaS04} on the use of
{\em relaxation and linearization techniques} \citep[cf.,][]{SherA94}, in combination with {\em spatial branching techniques} \citep[cf.,][]{TamaS04}, has lead to the development of
the two well-known global solvers {\tt Couenne} \citep{Belo10} and {\tt \BARON{}}~\citep{Sahi96} for NLPs. Another solver that combines these type of techniques, together with techniques to exploit the problems structure is {\tt GloMIQO}, developed by \cite{misener2013glomiqo} for the solution of more general quadratically constrained quadratic programs with integer variables.
More recently, specialized solution approaches have been developed for special classes of QP. In particular, \cite{bonami2016solving} develop a special branch-and-cut algorithm for box constrained QPs based on using cuts derived from the boolean quadric polytope. Also, \cite{bonami2016cut} develop new specialized cuts that are used within a spatial branch and bound algorithm  to solve standard QPs.
For further review of numerical and theoretical results on the solution of QPs using NLP techniques, we refer the reader to \citet{gould2003} and \citet{gao2004}.

Besides NLP techniques, convex optimization techniques \citep[cf.,][]{Rene01, ben-tal01} have also been used
to address the solution of~QPs.
For example, \citet{nesterov1998} and later \citet{kim2001, kim2003}, explored the use of semidefinite programming (SDP)
as well as second-order cone relaxations to approximately or globally solve a QP.

More recently,
\cite{burer2009} proposed a SDP-based branch and bound approach to globally solve box-constrained
QPs; they reformulate a QP by adding the QP's corresponding Karush-Kuhn-Tucker~(KKT) conditions as redundant constraints.
Let us refer to this quadratically constrained quadratic program (QCQP) as \KKTQCQP. To solve \KKTQCQP{}, \cite{burer2009} construct a finite {\em KKT-branching} tree by branching on the resulting problem's complementarity constraints. SDP relaxations of \KKTQCQP\ are used to obtain lower bounds at each node of the KKT-branching tree.
On the other hand, to obtain upper bounds a (local) QP-solver based on NLP techniques is used.
\cite{chen2012} improved the solution methodology of \cite{burer2009} by
obtaining tighter lower bounds at each node of the KKT-branching tree.
For that purpose,  the {\em double non-negative} (DNN) relaxation of the {\em completely positive} reformulation \citep{burer2009copo} of \KKTQCQP\ at each node of the KKT-branching tree is used.
\cite{chen2012} provide a {\tt MATLAB} implementation of their approach called {\tt {\tt quadprogBB}}. In this implementation, the {\tt MATLAB}  (local)  QP solver {\tt quadprog} is used to obtain the upper bounds while the algorithm proposed by \citet{BurerDNN} is used to obtain lower bounds, at each node of the KKT-branching tree. \citet{chen2012} show that this solution approach typically outperforms the solver {\tt Couenne} and the approach proposed by \cite{burer2009}
on a test bed of publicly available QP instances.
This makes the solver {\tt quadprogBB} a current benchmark for the global solution of QP problems.

The rest of the paper is organized as follows.
 In Section 2, we formally introduce the QP problem and present the theoretical
 results that serve as the foundation for the  proposed solution approach.
 In Section 3, we  illustrate the effectiveness of this approach by presenting relevant numerical
 results on test instances of the QP problem. To conclude, in Section 4, we provide conclusions and directions for future work.
\section{Solution Approach for non-convex QPs}
\label{sec:1}

We consider the following quadratic programming problem
\begin{equation}
\label{eq:QP}
\begin{array}{rrl}
{\rm QP: }\qquad\qquad\qquad&\min \quad& \frac{1}{2}x\tr Hx + f\tr x\\
&\st &Ax = b\\
&& x \ge 0,
\end{array}
\end{equation}
where $f\in \R^n$, $A\in \R^{m\times n}$, $b\in \R^m$, and $H \in \R^{n\times n}$ is a symmetric matrix. Note that there is no assumption on the matrix~$H$ being positive semidefinite; that is, QP is in general a non-convex optimization problem \citep[cf.,][]{Bert99}.

Similar to \citet{burer2009} and \citet{chen2012}, we assume that the feasible set of QP
is nonempty and bounded.
However, in what follows, no further assumption is made about the 
feasible set of QP.

\subsection{Mixed-integer linear programming reformulation}
\label{sec:SolApproach}
After introducing the Lagrange multipliers $\mu \in \R^m$ for its equality constraints
and  $\lambda \in \R^n$ for its non-negativity constraints,
the KKT conditions for QP are given by
\begin{subequations}
\label{eq:KKT}
\begin{align}
 Hx + f +  A\tr \mu - \lambda &= 0 \label{LKKT}\\
 x\tr  \lambda &=0 \label{complKKT}\\
Ax &=b \label{feasKKT}\\
 x \ge 0,  \lambda &\ge 0. \label{nonNegKKT}
\end{align}
\end{subequations}
In what follows, we will refer to the set
\begin{equation}
\label{KKTset}
\Lambda_{KKT}  = \{ (x,\mu,\lambda) \in \R^{2n+m}:  (x,\mu, \lambda) \text{ satisfy}~\eqref{LKKT}-\eqref{nonNegKKT} \}
\end{equation}
as the {\em KKT points} of QP.

Note that because the feasible set of QP~\eqref{eq:QP} is a polyhedron, the KKT conditions \eqref{eq:KKT} are first order necessary conditions for the optimal solutions of~QP~\citep[see, e.g.][Thm. 3.3]{eustaquio2008constraint}. Thus, one can add these KKT conditions as redundant constraints in QP
to obtain the following equivalent formulation of QP,
\begin{equation}
\label{QPr1}
\begin{array}{rll}
\min \quad& \frac{1}{2}x\tr Hx + f\tr x\\
\st &  Hx + f +  A\tr \mu - \lambda &= 0\\
&x\tr  \lambda =0 \\
&Ax = b \\
&x \ge 0,  \lambda \ge 0. \\
\end{array}
\end{equation}
As shown by \citet[][Thm. 2.4]{Gian73}, one can use the KKT conditions~\eqref{LKKT}--\eqref{feasKKT} to linearize the objective of~\eqref{QPr1}. Namely, for any feasible solution $x \in \R^n$ of~\eqref{QPr1}, we have
\begin{equation}
\frac{1}{2}x\tr Hx + f\tr x  =\frac{1}{2}(f\tr  x -x\tr  A\tr  \mu + x\tr \lambda) = \frac{1}{2}(f\tr  x -b\tr \mu).  \label{LOBJ}
\end{equation}
As a result, problem \eqref{QPr1} is equivalent to the following problem with a linear (instead of quadratic) objective.
\begin{equation}
\label{QPr2}
\begin{array}{rll}
\frac 12\min &  f\tr  x  - b\tr  \mu \\
\st &  Hx + f +  A\tr \mu - \lambda &= 0\\
&x\tr  \lambda =0 \\
&Ax = b \\
&x \ge 0, \lambda \ge 0. \\
\end{array}
\end{equation}
Notice that in \eqref{QPr2}, the complexity of QP is captured
in the complementary constraints~$x\tr  \lambda =0$.
 Next, we address the complementary  constraints in \eqref{QPr2} by using
 {\em {Big-M}} constraints. For that purpose, in Section~\ref{sec:primaldualbounds}, we derive upper bounds $U, V \in \R^n$
 on the decision variables~$x, \lambda \in \R^n$ of~\eqref{QPr2}
 such that there are (globally) optimal KKT points~$(x,\mu,\lambda) \in \R^{2n + m}$ of ~QP  satisfying $x \leq U$, $\lambda \leq V$.
 Using these upper bounds, one can show (see, Theorem~\ref{thm:boundDual}) that a global optimal solution of QP can be obtained by solving
 the following MILP
\begin{equation}
\label{eq:milpQP}
\begin{array}{rllll}
{\rm IQP: }\qquad\qquad\qquad\frac 12\min &  f\tr  x  - b\tr  \mu \\
\st &  Hx + f +  A\tr \mu - \lambda = 0\\
&Ax = b \\
&  0 \le x_j \le z_j U_j & j=1,\dots,m\\
&  0 \le \lambda_j \le (1-z_j) V_j & j=1,\dots,m\\
&  z_j \in \{0,1\}& j=1,\dots,m.
\end{array}
\end{equation}
Specifically, problem IQP is a MILP with the same optimal value as QP whose optimal solutions are optimal solutions of QP.



\subsection{Bounding the primal and dual variables}
\label{sec:primaldualbounds}

As mentioned earlier, the first step in obtaining problem IQP
is to derive explicit upper bounds $U, V \in \R^n$
 such that there are optimal KKT points $(x,\mu,\lambda) \in \R^{2n+m}$ of QP satisfying $x \leq U$, $\lambda \leq V$.

Similar to \citet{chen2012}, using the assumption that the feasible set of QP is non-empty and bounded, one can compute the upper bounds $U \in \R^n_+$ on the primal variables $x\in \R^n$ by setting:
\begin{equation}\label{eq:primalBound}
U_j := \max \{x_j : Ax = b,\, x\ge 0\},
\end{equation}
for every $j=1,\dots,n$.

Using assumptions stronger than ours, \citet{chen2012}  show that $\Lambda_{KKT}$, the set of KKT points,  is bounded.
As the following example illustrates, under our weaker assumptions, the set  $\Lambda_{KKT}$ could be unbounded.
\begin{example}\label{ex:dualunbdd}
{\em Consider the instance of  QP  defined by setting
\[H = { \begin{bmatrix} 2 & 0 & 0\\ 0 & -1 & 0\\ 0& 0& 1\end{bmatrix}}\qquad f = {\begin{bmatrix} 2 \\ 4\\ 3\end{bmatrix}}\qquad A = { \begin{bmatrix} 2 & 1 & 1\\ 1& 1& 1\end{bmatrix}} \qquad b = \begin{bmatrix} 1 \\ 1 \end{bmatrix}.\]
 Note that in this case the feasible region of  QP  is $\{[0,1-t,t]\tr :0\le t\le 1\}$, which is bounded and non-empty.
However, the set of KKT points $\Lambda_{KKT}$~\eqref{KKTset} is unbounded. Specifically, notice that for any $v \ge 1$ the following is a KKT point for  QP :
\[x =  {\small\begin{bmatrix} 0 \\ 1\\ 0\end{bmatrix}}\qquad \mu = {\small\begin{bmatrix} v\\ -3 -v \end{bmatrix}}\qquad \lambda = {\small\begin{bmatrix} -1+ v\\ 0\\ 0\end{bmatrix}}.\]}
\end{example}

Thus, to handle the complementarity contrains in \eqref{QPr2} using Big-M constraints, we do not try to obtain a bound for the value of the entries of $\lambda \in \R^m$ for all KKT points. Instead, in Theorem~\ref{thm:boundDual}, we prove that there exist a bound that we can impose in the dual variables, without discarding all (globally) optimal KKT points of  QP. For this purpose, we make use of fundamental  results on the approximate solution of systems of linear equations \cite[e.g.,][]{GuleHR95,  Mang81}.

Let us first define a particular instance of the well-known Hoffman bound~\citep{Hoff52}, closely following the notation in \cite{GuleHR95}.
\begin{definition}
\label{def:Hoffman}
Fix the norm $\|\cdot\|_{\alpha}$ on $\R^n$ and the norm $\|\cdot\|_{\beta}$ on $\R^m$. 
Given $A \in \R^{m \times n}$ and $b \in \R^m$, let $F := \{x \in \R^n_+: Ax = b\}$. Let $\hoffgen \in \R$ be the smallest constant
satisfying:
\begin{equation}
\label{eq:bound}
\text{For all } y \in \R^n \text{ such that }Ay = b, \text{ there is } x \in F \text{ such that } \|x - y\|_{\alpha}  \le \hoffgen \|y^-\|_{\beta}.
\end{equation}
\end{definition}

Above, for any $y \in \R^n$, $y^- \in \R^n$ is the vector difined by $y^-_i= \max\{0, -y_i\}$, $i=1,\dots,n$. That is, Definition~\ref{def:Hoffman}
corresponds to the Hoffman bound obtained when looking at perturbations of only the non-negative constraints of the polyhedron $F := \{x \in \R^n_+: Ax = b\}$.

In what follows we will use the following notation to denote the dual norm associated to a given norm.

\begin{definition}
{\em Given a norm $\|\cdot\|$ on $\R^n$, its associated {\em dual norm} on $\R^n$, denoted $\|\cdot\|^*$ is defined as:
\[
\|x\|^* = \sup\{ x\tr z: z\in \R^n,\,\|z\| \le 1\}.
\]}
\end{definition}
In particular, for any $x \in \R^n$, $\|x\|_{\infty}^* = \|x\|_1$.

Using Definition~\ref{def:Hoffman}, we provide in Theorem~\ref{thm:boundDual} below, the desired bound to be used on the dual variables of QP in the IQP formulation~\eqref{eq:milpQP}.

\begin{theorem}\label{thm:boundDual}
Let $A \in \R^{m \times n}$ and $b \in \R^m$ be such that the set  $F := \{x \in \R^n_+: Ax = b\}$ is non-empty and
bounded.
Let $\hoffgen$ be defined by  \eqref{eq:bound}, and  $M > (\kappa  + \|f\|_{\alpha}^* )\hoffgen \|e\|_{\beta}$, where $\kappa := \max \{\|Hx\|^*_{\alpha}: Ax = b, x\ge 0\}$.
Then, there exists an optimal KKT point $(x^*,\mu^*,\lambda^*)$ for~QP  such that
$e\tr \lambda^* \le M$,  and $\mu^* \in \R^m$.
\end{theorem}

\proof{Proof.}
Consider the following perturbed version of QP:
\begin{equation}
\label{perQP}
\begin{array}{ll}
\min \quad& \frac{1}{2}x\tr Hx + f\tr x + Mt\\
\st &Ax = b\\
& x \ge - te\\
& 0 \le t \le \delta,\\
\end{array}
\end{equation}
where $e$ is the vector of all ones, and $\delta > 0$.
Notice that the feasible set
of \eqref{perQP} is a closed subset of $\{x \in \R^n: Ax = b,\, x \ge -\delta e\}\times [0, \delta]$ which is non-empty and  bounded as $F$ is non-empty, bounded, and the recession cone of $\{x \in \R^n: Ax = b,\, x \ge -\delta e\}$ is equal to the recession cone of $F$. Thus, the optimal value of~\eqref{perQP} exists and it is attained.

Let $(x^*,t^*)$ be an optimal solution of~\eqref{perQP}. Then, there exists $(\mu^*, \lambda^*, \rho^*,\omega^*) \in \R^{n+m+2}$ such that $(x^*,t^*,\mu^*,\lambda^*,\rho^*,\omega^*)$ satisfies the KKT conditions associated with problem~\eqref{perQP}
\begin{equation}
\label{eq:KKT2}
\begin{array}{rllll}
 Hx^* + f +  A\tr \mu^* - \lambda^* &= 0\\
  M - e\tr \lambda^* -\rho^* + \omega^* &= 0\\
 (x^*-t^*e)\tr  \lambda^* &=0 \\
  {t^*} \rho^* &=0 \\
  (\delta-t^*) \omega^* &=0 \\
Ax^* &= b\\
 x^* + t^*e &\ge 0 \\
0 \le t^* &\le \delta\\
 \lambda^*,\rho^*,\omega^* &\ge 0,
\end{array}
\end{equation}
where $\mu^* \in \R^m$, $\lambda^* \in \R^n$,  $\rho^* \in \R$, $\omega^* \in \R$, are respectively the Lagrangian multipliers of problem~\eqref{perQP} associated
with the linear constraints, lower bounds in the decision variables~$x \in \R^n$, and lower and upper bounds on the decision variable $t \in \R$.

Now we claim that $t^* = 0$. In that case, notice that the complementarity
constraint $ (\delta-t^*) \omega^* =0$ in \eqref{eq:KKT2} implies $\omega^* = 0$ and thus,
from the equation $M - e\tr \lambda^* -\rho^* + \omega^* = 0$ in~\eqref{eq:KKT2} and the fact that $\rho^* \ge 0$, it follows that $(x^*,\mu^*,\lambda^*)$ satisfies the statement of the theorem.

To show that $t^* = 0$, note that from Definition~\ref{def:Hoffman}, it follows that there exists $x' \in F$ such that\begin{equation}
\label{eq:genbound}
\|x'-x^*\|_{\alpha} \le \hoffgen t^*\|e\|_{\beta}.
\end{equation}
In problem~\eqref{perQP}, $(x',0)$ is a
 feasible solution and thus the objective value of $(x',0)$  is no smaller than the objective value of~$(x^*,t^*)$. That is,
\begin{equation}\label{eq:optxu}
\frac{1}{2}{x^*}\tr Hx^* + f\tr x^* + Mt^* \le \frac{1}{2}{x'}\tr Hx' + f\tr x'.
\end{equation}
Therefore,
\begin{align*}
Mt^*&\le \frac{1}{2}({x'}\tr Hx'-{x^*}\tr Hx^*) + f\tr (x'-x^*)                                    \\
&=\frac{1}{2}(x'+x^*)\tr H(x'-x^*) + f\tr (x'-x^*) \\
& \le \frac 12 \|H(x'+x^*)\|^*_{\alpha}\|x'-x^*\|_{\alpha} +   \|f\|_{\alpha}^*\|x'-x^*\|_\alpha \\
& \le \left ( \frac{1}{2}( \|Hx'\|^*_{\alpha} + \|Hx^*\|^*_{\alpha})+ \|f\|_{\alpha}^* \right ) \|x'-x^*\|_\alpha.
\end{align*}
Thus, using~\eqref{eq:genbound} we have
\begin{equation}
\label{eq:last}
Mt^* \le \left ( \frac{1}{2}( \kappa + \kappa_{\delta})+ \|f\|_{\alpha}^* \right ) \hoffgen \|e\|_{\beta} t^*
\end{equation}
where $
\kappa_{\delta} := \max \{\|Hx\|^*_{\alpha}: Ax = b, x \ge -\delta e, x \in \R^n\}$.

As $M > \left (\kappa + \|f\|_{\alpha}^* \right ) \hoffgen \|e\|_{\beta}$, and $\kappa_{\delta} \downarrow \kappa$ when $\delta \downarrow 0$, taking $\delta>0$ small enough we have that $M > \left ( \frac{1}{2}( \kappa + \kappa_{\delta})+ \|f\|_{\alpha}^* \right ) \hoffgen \|e\|_{\beta}$. Thus,~\eqref{eq:last}  implies that $t^* = 0$.
\endproof{}

\subsection{Computation of the dual bounds}
\label{sec:comphauff}

Theorem~\ref{thm:boundDual} provides the bounds needed for the reformulation of QP as IQP~\eqref{eq:milpQP} in terms of the Hoffman constant $\hoffgen$ introduced in Definition~\ref{def:Hoffman}. Next, we discuss how this constant can be obtained in closed-form for important special classes of QP, as well as how it can be computed for general classes of QP. 

\subsubsection{Standard Quadratic Programming.}
\label{sec:SQPbound}
Consider the {\em standard quadratic program}~(SQP):
\begin{equation}
\label{eq:SQP}
\begin{array}{llllll}
{\rm SQP:}\quad\quad\quad\dmax_{x \in \Delta}& \frac{1}{2}x\tr Hx + f\tr x\\
\end{array}
\end{equation}
where
$\Delta = \left \{x\in\mathbb{R}^n:\sum_{i = 1}^n x_i = 1, x\geq 0 \right \}$,
 is the standard simplex. The SQP problem is fundamental in optimization and arises in
many applications \citep[see, e.g.,][]{Bomze98}.
Next,  we show that in this case, the Hoffman bound $\hoffgen$ introduced in Definition~\ref{def:Hoffman} can be computed in closed-form for suitable choices of the norm $\| \cdot \|_{\alpha}$ on $\R^n$ and the norm $\| \cdot \|_{\beta}$ on $\R^m$.

\begin{proposition}
\label{prop:sqpbound} 
Consider the norm $\|\cdot\|_{\alpha} = \|\cdot\|_{1}$ on $\R^n$ and the norm $\|\cdot\|_{\beta} = \|\cdot\|_{1}$ on~$\R^m$.
Let $A = e\tr $ and  $b=1$. Then $\hoff = 2$.
\end{proposition}

\proof{Proof.}
Let  $y \in \R^n$ such that $e\tr y =1$ be given.
Let $I= \{i \in \{1,\dots,n\}: y_i \ge 0\}$ and $I^c = \{1,\dots,n\} \setminus I$, and consider the case $I^c \neq \emptyset$ (otherwise, the statement follows by letting $x = y$ in Definition~\ref{def:Hoffman}). Note that $e\tr y = 1$ implies $I \not = \emptyset$ and
that $\sum_{i \in I} y_i  = 1 -  \sum_{i \in I^c} y_i = 1 + \|y^-\|_1$.
Let $x \in \R^n$ be defined by setting $x_i = 0$ for all $i \in I^c$, and $x_i = \frac{1}{1 + \|y^-\|_1}y_i$ for all~$i \in I$. Clearly, $x \in \Delta$. Furthermore, for any $i \in I^c$, $|x_i - y_i| = -y_i$.
Also for any $i \in I$, we have $|x_i - y_i| = \frac{\|y^-\|_1}{1 + \|y^-\|_1}y_i$.
Thus, $\|x-y\|_{1} = -\sum_{i \in I^c} y_i + \frac{\|y^-\|_1}{1 + \|y^-\|_1}\sum_{i \in I} y_i = 2\|y^-\|_1$. That is, $\hoff \le 2$. To show $\hoff \ge 2$, consider $y = (n, -1, \dots, -1)$.
For any $x \in \Delta$, it follows that $\|x-y\|_{1} = |x_1-n| + \sum_{i=2}^n |x_i+1|  = 2n-1 -x_1 + \sum_{i=2}^n x_i = 2(n-x_1) \ge 2(n-1) = 2\|y^-\|_1$. 
 \endproof{}

 Thus, \eqref{eq:SQP} can be reformulated as IQP by letting:
\begin{equation}
\label{eq:sqpbd}
U = e \text{ and } V>  Me,
\end{equation}
with
\begin{equation}
\label{eq:sqpbdM}
M =  2n \left(\|H\|_{\infty,\infty}+ \|f\|_{\infty} \right ),
\end{equation}
where we have used that
\[\kappa = \max\{\|Hx\|_{\infty}: e\tr x = 1, x\in \R^n_+\} \le  \max_{i,j \in \{1,\dots,n\}} |H_{ij}| =: \|H\|_{\infty,\infty}.\]

\begin{remark}
It is worth mentioning that a proof similar to the one given in Proposition~\ref{prop:sqpbound} shows that
if~$\|\cdot\|_1$  is replaced with $\|\cdot\|_\infty$ in Proposition~\ref{prop:sqpbound}, the corresponding Hoffman constant would be equal to $n-1$. However, this leads to a weaker bound $V$ than the one given in~\eqref{eq:sqpbd}.
\end{remark}

\subsubsection{Quadratic Programming with Box Constraints.}
\label{sec:boxqp}
Now, consider the  {\em box-constrained QP} (BoxQP) 
\begin{equation}
\label{eq:BQP}
\begin{array}{llllll}
\rm{BoxQP:}\qquad\qquad\qquad&\max\qquad & \frac{1}{2}x\tr Hx + f\tr x\\
&\st& l \le x \le  u,
\end{array}
\end{equation}
where  $l, u \in \R^n$ are given bounds on the primal variables of BoxQP satisfying (w.l.o.g.) $l < u$ (component-wise).
Problem BoxQP is equivalent to the following   QP  problem:
\begin{equation}
\label{eq:BQP_std}
\begin{array}{llllll}
\max & \frac{1}{2}x\tr Hx + (Hl+f)\tr x\\
\st& x + s =  u - l\\
& x\ge 0 , s \ge 0.
\end{array}
\end{equation}

Next,  we show that in this case, the
Hoffman bound $\hoffgen$ introduced in Definition~\ref{def:Hoffman} can be computed in closed-form for a suitable choice 
 of the norm $\| \cdot \|_{\alpha}$ on $\R^n$ and the norm $\| \cdot \|_{\beta}$ on $\R^m$.

\begin{proposition}
\label{prop:box}
Consider the norm $\|\cdot\|_{\alpha} = \|\cdot\|_{\infty}$ on $\R^n$ and the norm $\|\cdot\|_{\beta} = \|\cdot\|_{\infty}$ on~$\R^m$.
Let $I$ denote the identity matrix in~$\R^{n\times n}$,
$b \in \R^n_{+}$, and $A = [I, I]$. Then $\hoff = 1$.
\end{proposition}

\proof{Proof.}
Let $(y, z) \in \R^{2n}$ such that $y + z = b$ be given. Define $x = y^+ -z^-$ and $s = z^+ - y^-$. We claim $(x,s) \in F = \{ (x,s) \in \R^{2n}_+: x+ s = b\}$. To show this notice first that $x+s = y^+ -z^- + z^+ - y^- = y + z = b$. Now let $i \in \{1,\dots,n\}$. If $z_i^- = 0$ then $x_i = y^+_i \ge 0$. Thus assume $z_i^- > 0$. Then $z_i^+ = 0$ and  $x_i = b_i - s_i = b_i + y_i^- \ge 0$. Thus $x \ge 0$. Similarly $s \ge 0$.
To finish, notice that $\|(y, z) - (x,s)\|_\infty = \|(-y^- +z^-, -z^- +y^- )\|_\infty = \|(y^-, z^-)\|_\infty = \|(y, z)^-\|_\infty$.
This shows that $\hoff \le 1$. To show $\hoff \ge 1$, let $y = -e$ and $z = b+e$.
For any $(x,s) \in F$, it follows that $\|(x,s)-(y,z)\|_{\infty} \ge |x_1+1| \ge 1 = \|(y,z)^-\|_{\infty}$. 
\endproof{}

Using Proposition \ref{prop:box}, we obtain that \eqref{eq:BQP_std} can be reformulated as IQP by letting:
\begin{equation}
\label{eq:boxbd}
U = {\small \begin{bmatrix} u-l \\ u-l\end{bmatrix}} \text{ and } V > Me,
\end{equation}
with
\begin{equation}
\label{eq:boxbdM}
M = \min(n\|H\|_{\infty,\infty}\|u-l\|_1, \|H\|_{1,1}\|u-l\|_\infty) + \|f +Hl\|_1,
\end{equation}
where we have used that
\[\kappa = \max\{\|Hx\|_{1}: x+s = u-l, x,s\in \R^n_+\} \le \min(n\|H\|_{\infty,\infty}\|u-l\|_1, \|H\|_{1,1}\|u-l\|_\infty),\]
where $\|H\|_{1,1}:=\sum_{i \neq j \in \{1, \dots, n\}}|H_{ij}|$.

\subsubsection{General  Quadratic Programming.}
\label{sec:genQPbound}
Note that to compute an appropriate $M$ value in Theorem~\ref{thm:boundDual}, it is enough to let
$M > (\kappa  + \|f\|_{\alpha}^* )\hoffpgen \|e\|_{\beta}$ for some constant $\hoffpgen \ge \hoffgen$. As shown bellow, $\hoffpgen$ can be computed in general using \citet[][Theorem 3.2]{GuleHR95}.

\begin{proposition}
\label{prop:Hoffman}
Fix the norm $\|\cdot\|_{\alpha}$ on $\R^n$ and the norm $\|\cdot\|_{\beta}$ on $\R^m$. 
Let $A \in \R^{m \times n}$ and $b \in \R^m$ be such that the set  $F := \{x \in \R^n_+: Ax = b\} \not = \emptyset$. Also, let
\begin{equation}
\label{eq:boundset}
\bar{\sigma}(A) = \left \{ (\mu^+, \mu^-, \lambda) \in \R^{m+n}: \|A\tr (\mu^+-\mu^-) - \lambda\|_{\alpha}^* \le 1, \mu^+,\mu^- \in \R_+^m, \lambda \in \R^n_+ \right \},
\end{equation}
and
\begin{equation}
\label{eq:boundpena}
\hoffpgen = \max\{\|(\mu^+,\mu^-, \lambda)\|^*_{\beta} : (\mu^+,\mu^-, \lambda) \text{ is an extreme point of } \bar{\sigma}(A)\}.
\end{equation}
Then, $\hoffpgen \ge \hoffgen$.
\end{proposition}

\proof{Proof.}
First notice that for all $y \in \R^n$,
\[
\min_{x \in F} \| x - y\|_{\alpha} = \min \{\| x - y\|_{\alpha}: A'x \le b', x \in \R^n\}
\]
where
\[
A' = \begin{bmatrix} A \\ -A \\ -I \end{bmatrix}, b'= \begin{bmatrix} b \\ -b \\ 0 \end{bmatrix}.
\]
Thus, it follows from \citet[][Theorem 3.2]{GuleHR95} that
\begin{equation}
\label{eq:guler}
\min_{x \in F} \| x - y\|_{\alpha} \le \hoffpgen \left \| \begin{pmatrix} (Ay - b)^+\\ (Ay - b)^- \\ y^- \end{pmatrix} \right \|_{\beta},
\end{equation}
after identifying $\hoffpgen = K_{\alpha\beta}(A')$, $\bar{\sigma}_{\alpha}(A) = \sigma_{\alpha}(A')$~\citep[see,][Theorem 3.2]{GuleHR95}.

From \eqref{eq:guler}, it follows that for any $y \in \R^n$ such that $Ay = b$, then
\begin{equation}
\label{eq:ourguler}
\min_{x \in F} \| x - y\|_{\alpha} \le \hoffpgen \left \| \begin{pmatrix} 0\\ 0 \\ y^- \end{pmatrix} \right \|_{\beta}  = \hoffpgen \left \|  y^-  \right \|_{\beta}.
\end{equation}
Also, from \eqref{eq:ourguler}, $\hoffpgen \ge \hoffgen$ follows from Definition~\ref{def:Hoffman}, as $\hoffgen$ is the smallest constant satisfying~\eqref{eq:ourguler}.
\endproof{}

In order to use Proposition~\ref{prop:Hoffman} to reformulate QP~\eqref{eq:QP} as IQP~\eqref{eq:milpQP}, one can first, similar to~\cite{chen2012}, normalize the primal variables of QP to be between $0$ and $1$. Namely, under the boundedness assumption considered here, on has that QP is equivalent to:
\begin{equation}
\label{eq:QPnormalized}
\begin{array}{rrl}
\qquad\qquad\qquad&\min \quad& \frac{1}{2}x\tr \tilde{H}x + \tilde{f}\tr x\\
&\st &\tilde{A}x = b\\
&& x \ge 0,
\end{array}
\end{equation}
where $\tilde{A}_{ij} := A_{ij}U_{i}U_{j}$, $\tilde{f}_i := f_iU_i$, and $\tilde{A}_{ij} := A_{ij}U_j$ for all $i,j \in\{1,\dots,n\}$, and $U \in \R^n_+$ is given by~\eqref{eq:primalBound}. Now, using Proposition~\ref{prop:Hoffman}, 
and choosing the norm $\|\cdot\|_{\alpha} = \|\cdot\|_{\infty}$	on $\R^n$ and the norm $\|\cdot\|_{\beta} = \|\cdot\|_{\infty}$	on $\R^m$,
we obtain that QP~\eqref{eq:QP} can be reformulated as IQP~\eqref{eq:milpQP} by letting $H = \tilde{H}$, $f = \tilde{f}$, $A = \tilde{A}$, 
\begin{equation}
\label{eq:genQPbd}
U = e, \text{and } V >  Me,
\end{equation}
with
\begin{equation}
\label{eq:genQPbdM}
M=  \left (\|\tilde{H}\|_{1,1} + \|\tilde{f}\|_{1} \right )\hoffpgen,
\end{equation}
where we have used that $\kappa = \max\{\|\tilde{H}x\|_{1}: x \le e, x \in \R^n_+\} \le \|\tilde{H}\|_{1,1}$.

In the case that one chooses  the norm $\|\cdot\|_{\alpha} = \|\cdot\|_{1}$	on $\R^n$ and the norm $\|\cdot\|_{\beta} = \|\cdot\|_{1}$	on~$\R^m$, then $M =  n(n\|\tilde{H}\|_{\infty,\infty} + \|\tilde{f}\|_{\infty})\hoffpgen$. However, empirical results on test instances shows that this latter $M$ is weaker than the bound obtained using~\eqref{eq:genQPbdM}.


Obtaining an efficient way to compute the constant $\hoffpgen$ in~\eqref{eq:boundpena} is however still an open question \citep[see, e.g.,][]{Ng04, GuleHR95, Penahoff17}. For illustrative purposes, in Section~\ref{sec:compbounds}, we show the results of using an algorithm recently proposed in~\citet{Penahoff17} to compute $\hoffpgen$, and the corresponding bound $M$ in~\eqref{eq:genQPbdM}.

An alternative and efficient way to compute bounds on the dual variables of a general instance of QP, is to use the bounds on the dual variables proposed by \citet[][Proposition~3.1]{chen2012} which are valid for  QP  instances having a strictly non-negative feasible solution (i.e., a feasible solution satisfying $x >0$), and
can be computed by solving a LP \citep[cf.,][eq. (19)]{chen2012}. Specifically, notice that after obtaining the primal bounds $U\in \R^n$ using \eqref{eq:primalBound},
problem  QP  is equivalent to
\begin{equation}
\label{eq:boundedQP}
\begin{array}{rl}
\min \quad& \frac{1}{2}x\tr Hx + f\tr x\\
\st &Ax = b\\
& 0 \le x \le U.\\
\end{array}
\end{equation}
Following the notation used thus far and letting $\rho \in \R^n$ be the dual variables associated with the upper bound constraints on the
variables $x \in \R^n$ in \eqref{eq:boundedQP}, it follows from the KKT conditions of \eqref{eq:boundedQP} that any of its optimal solutions must satisfy:
\begin{subequations}
\label{eq:boundedKKT}
\begin{align}
 Hx + f +  A\tr \mu - \lambda + \rho &= 0 \label{bdLKKT}\\
 x\tr  \lambda =0, (U-x)\tr \rho & = 0 \label{bdcomplKKT}\\
 Ax &=b \label{bdfeasible} \\
\notag x \ge 0, \lambda \ge 0, \rho &\ge 0.
\end{align}
\end{subequations}
Also, after multiplying \eqref{bdLKKT} by a feasible solution $x \in \R^n$ of \eqref{eq:boundedQP} and using \eqref{bdcomplKKT},
\eqref{bdfeasible}, it follows that any optimal solution of
\eqref{eq:boundedQP} also satisfies:
\begin{equation}
x\tr Hx + f\tr x +  b\tr \mu  + U\tr \rho = 0.
\end{equation}
Then, if  QP  has a feasible solution $x \in \R^n$ satisfying $x_i > 0$, $i=1,\dots,m$, it
follows from \citet[][Proposition~3.1]{chen2012} that bounds on the dual variables $V \in \R^n$ required for the MILP reformulation IQP of~ QP
can be computed by solving the following LP:
\begin{equation}
\label{eq:genbd}
V_j = \max \left \{\lambda_j: \begin{array}{l}
                                           Hx + f +  A\tr \mu - \lambda + \rho = 0 \\
                                           H \bullet X + f\tr x +  b\tr \mu  + U\tr \rho = 0\\
                                            0 \le X_{ij} \le U_iU_j, i,j=1,\dots,n\\
                                           0 \le x \le U, \lambda \ge 0, \rho \ge 0, X \in \S^n \\
                                           \end{array} \right \},
\end{equation}
where $H \bullet X$ indicates the trace of the matrix $HX$, the matrix $X \in S^n$ represents the {\em linearization}
of the matrix $xx\tr  \in \S^n$, and $\S^n$ is the set of $n \times n$ real symmetric matrices.

\begin{remark}
{\em In \citet[][]{chen2012}, eq. (18) is used to refine the dual variable
bounds after scaling the problem so that its variables are between zero and one. However, this refinement
of the bounds is not necessary to obtain their result in Proposition~3.1. The refined version of these bounds is however the one implemented in \qIP, the implementation of our solution approach for general instances of QP.}
\end{remark}

\begin{remark}{\em
From Theorem~\ref{thm:boundDual}, the constraint $e\tr \lambda \le M$ could be added in the IQP formulation of QP~\eqref{eq:milpQP}. Adding this constraint however has not led to improved solution times of IQP. Thus, this constraint is not used in the implementation of our proposed solution approach for QP.}
\end{remark}

\subsubsection{Bound comparison}
\label{sec:compbounds}

In light of the bounds on the dual variables discussed here and the dual bounds proposed by~\citet[][eq. (19)]{chen2012}, it is natural to compare their values and computing times for different instances of QP. Before doing so, however, it is important to emphasize that the dual bounds proposed here can be imposed on QP, even if the dual feasible set of QP is unbounded.

\begin{example}[Example~\ref{ex:dualunbdd} revisited]
\label{ex:dualnbdd2}{\em 
Recall the problem discussed in Example~\ref{ex:dualunbdd}, and consider the norm $\|\cdot\|_{\alpha} = \|\cdot\|_{\infty}$ on $\R^n$ and the norm $\|\cdot\|_{\beta} = \|\cdot\|_{\infty}$ on~$\R^m$. It is not difficult to see that in this case $\hoffgen \le 1$, and $\kappa = 1$. Since $\|f\|_1 = 9$, and $\|e\|_{\infty} = 1$, it follows that $M=(\kappa + \|f\|_1)\hoffgen \le 10$. Thus, for this problem we can bound the dual variables with the constraint 
\begin{equation}
\label{eq:exbound}
[ \lambda_1, \lambda_2, \lambda_3 ]\tr \le 11[ 1, 1,1 ]\tr.
\end{equation}
In fact, notice that the following optimal solution of the problem 
\[x^* =  {\small\begin{bmatrix} 0 \\ 1\\ 0\end{bmatrix}}\qquad \mu^* = {\small\begin{bmatrix} 0\\ -3  \end{bmatrix}}\qquad \lambda^* = {\small\begin{bmatrix} 0\\ 0\\ 0\end{bmatrix}},\]
satisfies the dual bounds~\eqref{eq:exbound}. On the other hand, from Example~\ref{ex:dualunbdd}, it follows that 
\[
\max\{\lambda_1: (x,\mu,\lambda) \in \Lambda_{KKT}\} = \infty.
\]
Thus, dual bounds for this problem cannot be computed using~\citet[][eq.~(19)]{chen2012}.
}
\end{example}

In Table~\ref{tab:compSQP}, we compare the bounds obtained in Section~\ref{sec:SQPbound} with the bounds obtained using~\citet[][eq. (19)]{chen2012} for a number of randomly selected SQP instances. From Table~\ref{tab:compSQP}, it is clear that the bounds obtained using Theorem~\ref{thm:boundDual}, and specifically, eq.~\eqref{eq:sqpbdM} are tighter than the ones obtained using~\citet[][eq. (19)]{chen2012} (and labeled ``RLT bounds'' in Table~\ref{tab:compSQP} for the {\em reformulation linearization techniques} (RLT) used to derive them) on a random sample of SPQ instances. In fact, this is the case for all the SQP instances considered in Section~\ref{sec:3}. 

\begin{table}[!htb]
\begin{center}
{\small
\begin{tabular}{lrrcrr}
\toprule
   & \multicolumn{2}{c}{Thm.~\ref{thm:boundDual}  bound} && \multicolumn{2}{c}{RLT bound} \\
   \cline{2-3} \cline{5-6}
\multicolumn{1}{c}{SQP Instance} & \multicolumn{1}{c}{Value} & \multicolumn{1}{c}{Time (s)} && \multicolumn{1}{c}{Value} & \multicolumn{1}{c}{Time (s)}\\
\midrule
spar030-060-1.mat:	&	 5,520 	&	0.0000	&&	 10,628 	&	0.4952	\\
spar030-070-3.mat:	&	 5,880 	&	0.0000	&&	 16,448 	&	0.4694	\\
spar050-040-1.mat:	&	 9,200 	&	0.0000	&&	 18,925 	&	1.4016	\\
spar050-050-3.mat:	&	 9,800 	&	0.0000	&&	 29,719 	&	1.8674	\\
spar060-020-1.mat:	&	 11,040 	&	0.0000	&&	 13,255 	&	1.8344	\\
spar070-075-1.mat:	&	 13,440 	&	0.0000	&&	 90,316 	&	17.0129	\\
spar080-025-2.mat:	&	 15,680 	&	0.0001	&&	 31,792 	&	13.8083	\\
spar080-025-3.mat:	&	 15,040 	&	0.0001	&&	 38,115 	&	13.8468	\\
spar090-025-2.mat:	&	 17,640 	&	0.0001	&&	 44,646 	&	22.8278	\\
spar090-025-3.mat:	&	 16,920 	&	0.0001	&&	 49,459 	&	23.3063	\\
spar090-050-3.mat:	&	 17,460 	&	0.0001	&&	 95,726 	&	37.7146	\\
spar090-075-1.mat:	&	 17,280 	&	0.0001	&&	 148,416 	&	49.4807	\\
spar100-050-2.mat:	&	 19,600 	&	0.0001	&&	 115,192 	&	67.8767	\\
spar100-050-3.mat:	&	 19,400 	&	0.0001	&&	 119,634 	&	67.6607	\\
\bottomrule
\end{tabular}
}
\end{center}
\caption{Comparison of bounds on $e\tr \lambda$ obtained using~\eqref{eq:sqpbdM} (	``Thm.~\ref{thm:boundDual}  bound'' columns) vs. using ~\citet[][eq. (19)]{chen2012} (``RLT bound'' columns), together with corresponding computation times for a random sample of SQP instances.}
\label{tab:compSQP}
\end{table}

Similarly, in Table~\ref{tab:compBoxQP}, we compare the bounds obtained in Section~\ref{sec:boxqp} with the bounds obtained using~\citet[][eq. (19)]{chen2012} for a number of randomly selected SQP instances. From Table~\ref{tab:compBoxQP}, it is clear that the bounds obtained using Theorem~\ref{thm:boundDual}, and specifically, eq.~\eqref{eq:boxbdM} are tighter than the ones obtained using~\citet[][eq. (19)]{chen2012} on a random sample of BoxPQ instances. In fact, this is the case for all the BoxQP instances considered in Section~\ref{sec:3}. Both in Table~\ref{tab:compSQP} and Table~\ref{tab:compBoxQP} the time differences are a result of the bounds resulting from Theorem~\ref{thm:boundDual} being computed from the closed-form formulas~\eqref{eq:sqpbdM},~\eqref{eq:boxbdM}, while the RLP bounds are obtained by solving a linear program~\citep[][eq. (19)]{chen2012}.

\begin{table}[!htb]
\begin{center}
{\small
\begin{tabular}{lrrcrr}
\toprule
   & \multicolumn{2}{c}{Thm.~\ref{thm:boundDual}  bound} && \multicolumn{2}{c}{RLT bound} \\
   \cline{2-3} \cline{5-6}
\multicolumn{1}{c}{BoxQP Instance} & \multicolumn{1}{c}{Value} & \multicolumn{1}{c}{Time (s)} && \multicolumn{1}{c}{Value} & \multicolumn{1}{c}{Time (s)}\\
\midrule
spar020-100-3.mat:	&	 9,708 	&	0.0001	&	&	 73,497 	&	0.2955	\\
spar030-060-1.mat:	&	 13,097 	&	0.0001	&	&	 152,505 	&	0.5289	\\
spar030-070-1.mat:	&	 14,887 	&	0.0002	&	&	 172,151 	&	0.5001	\\
spar030-070-2.mat:	&	 15,909 	&	0.0001	&	&	 176,093 	&	0.5362	\\
spar030-070-3.mat:	&	 16,827 	&	0.0001	&	&	 184,397 	&	0.4815	\\
spar030-080-1.mat:	&	 18,259 	&	0.0001	&	&	 219,338 	&	0.4578	\\
spar030-080-2.mat:	&	 18,532 	&	0.0001	&	&	 205,091 	&	0.4539	\\
spar030-080-3.mat:	&	 18,585 	&	0.0001	&	&	 202,502 	&	0.5233	\\
spar040-060-3.mat:	&	 25,889 	&	0.0001	&	&	 388,914 	&	0.7893	\\
spar040-080-1.mat:	&	 31,929 	&	0.0001	&	&	 524,796 	&	0.9596	\\
spar040-090-2.mat:	&	 37,109 	&	0.0001	&	&	 589,155 	&	1.1884	\\
spar070-025-1.mat:	&	 30,162 	&	0.0182	&	&	 819,461 	&	7.1339	\\
\bottomrule
\end{tabular}
}
\end{center}
\caption{Comparison of bounds on $e\tr \lambda$ obtained using~\eqref{eq:boxbdM} (``Thm.~\ref{thm:boundDual}  bound'' columns) vs. using~\citet[][eq. (19)]{chen2012} (``RLT bound'' columns), together with corresponding computation times for a random sample of BoxQP instances.}
\label{tab:compBoxQP}
\end{table}

In Table~\ref{tab:compgenQP}, we compare the bounds obtained in Section~\ref{sec:genQPbound} with the bounds obtained using~\citet[][eq. (19)]{chen2012} for a number of randomly selected general QP instances. From Table~\ref{tab:compgenQP}, it is clear that the bounds obtained using Theorem~\ref{thm:boundDual}, and specifically, eq.~\eqref{eq:genQPbdM} are weaker than the ones obtained using~\citet[][eq. (19)]{chen2012} on a random sample of general QP instances. In fact, this is the case for all the general QP instances considered in Section~\ref{sec:3}. 
In Table~\ref{tab:compgenQP} the time differences are a result of the bounds resulting from Theorem~\ref{thm:boundDual} being computed using an algorithm whose complexity is exponential on the size of the constraint matrix of the problem~\citep{Penahoff17}, while the RLP bounds are obtained by solving a linear program~\citep[][eq. (19)]{chen2012}.

\begin{table}[!htb]
\begin{center}
{\small
\begin{tabular}{lrrcrr}
\toprule
   & \multicolumn{2}{c}{Thm.~\ref{thm:boundDual}  bound} && \multicolumn{2}{c}{RLT bound} \\
   \cline{2-3} \cline{5-6}
\multicolumn{1}{c}{General QP Instance} & \multicolumn{1}{c}{Value} & \multicolumn{1}{c}{Time (s)} && \multicolumn{1}{c}{Value} & \multicolumn{1}{c}{Time (s)}\\
\midrule
st\_e26.mat	&	 49,200 	&	0.0575	&	&	 8,828 	&	0.0742	\\
st\_fp4.mat	&	 830,297 	&	1.6961	&	&	 594 	&	0.1380	\\
st\_fp5.mat	&	 47,263,557 	&	125.8385	&	&	 792 	&	0.2807	\\
st\_glmp\_kky.mat	&	 83,410 	&	22.4100	&	&	 339 	&	0.1013	\\
st\_glmp\_ss1.mat	&	 33,757 	&	6.5787	&	&	 429 	&	0.1067	\\
st\_m1.mat	&	 556,912,094 	&	0.0655	&	&	 19,060,333 	&	0.3594	\\
st\_pan2.mat	&	 6,017 	&	0.0544	&	&	 1,494 	&	0.0939	\\
st\_jcbpaf2.mat:	&	 \multicolumn{1}{r}{-}	&	\multicolumn{1}{r}{-}	&	&	 96,969 	&	0.2757	\\
st\_ph10.mat	&	 1,320 	&	0.5757	&	&	 27 	&	0.0679	\\
st\_ph2.mat	&	 69,951 	&	0.0601	&	&	 8,043 	&	0.1132	\\
st\_qpc\_m0.mat	&	 372 	&	0.0573	&	&	 35 	&	0.0501	\\
qp20\_10\_2\_1.mat	&	 69,846 	&	0.0671	&	&	 2,500 	&	0.5206	\\
qp30\_15\_1\_4.mat	&	 44,451 	&	0.0621	&	&	 1,661 	&	0.7943	\\
qp30\_15\_2\_4.mat	&	 41,625 	&	0.0602	&	&	 1,122 	&	0.9512	\\
qp40\_20\_1\_4.mat	&	 85,008 	&	0.0548	&	&	 8,700 	&	1.3852	\\
qp40\_20\_4\_1.mat	&	 523,052 	&	0.0606	&	&	 16,371 	&	3.3573	\\
qp50\_25\_1\_2.mat	&	 149,684 	&	0.0668	&	&	 12,476 	&	2.4781	\\
qp50\_25\_1\_4.mat	&	 169,868 	&	0.0676	&	&	 17,623 	&	2.1617	\\
\bottomrule
\end{tabular}
}
\end{center}
\caption{Comparison of bounds on $e\tr \lambda$ obtained using~\eqref{eq:genQPbdM} (``Thm.~\ref{thm:boundDual}  bound'' columns) vs. using~\citet[][eq. (19)]{chen2012} (``RLT bound'' columns), together with corresponding computation times for a random sample of (general) QP instances. Dash ``-'' indicates that the time limit of 1800 sec has been reached without computing the bound.}
\label{tab:compgenQP}
\end{table}

From the results in Table~\ref{tab:compSQP}, Table~\ref{tab:compBoxQP}, and Table~\ref{tab:compgenQP}, it is clear that using the bound of Theorem~\ref{thm:boundDual} can lead to tighter bounds on the QP dual variables $\lambda \in \R^n_+$ than the ones obtained using~\citet[][eq. (19)]{chen2012} when a tight bound on the Hoffman constant $\hoffgen$ used in Theorem~\ref{thm:boundDual} can be computed efficiently.

As illustrated in Example~\ref{ex:dualnbdd2}, the dual QP bounds obtained from Theorem~\ref{thm:boundDual} can be used even if the dual feasible set of QP is unbounded. In such case, it is not possible to use the \qBB solution methodology proposed by \citet{chen2012} to solve the problem, as the methodology requires (through a condition on the primal QP problem) the dual feasible set of QP to be bounded. To illustrate this (see Table~\ref{tab:unbounded}), we modify some general QP test instances in a simple way to make their dual feasible set unbounded. The modification we use is to pick the first variable $x_1$ of the instance and add the constraint $x_1 = x_1^*$, where $x_1^*$ is the value of $x_1$ in an optimal solution of the problem
 (i.e., this results in a problem that likely violates the interior condition required by \citet[see,][preceding Prop. 3.1]{chen2012}). As shown in Table~\ref{tab:unbounded}, these modified instances can be correctly solved using the approach proposed here with the bounds~\eqref{eq:genQPbdM}, while \qBB of  \citet{chen2012} is unable to solve them due to the unboundedness of some of the dual variables of the modified version of the problem. Specifically, Table~\ref{tab:unbounded} provides the name of the original instance (1st column), its optimal value (2nd column), the constraint that is added to the problem to make its dual feasible set unbounded while leaving its optimal value unchanged (3rd column), the value of the $M$ bound~\eqref{eq:genQPbdM} computed as a bound for the dual variables while still retaining at least an optimal solution (4th column), the optimal solution for the modified version of the instance obtained with \qIP (5th column), and the number of the dual variable of the modified version of the instance that \qBB detects to be unbounded which results in \qBB not being able to solve the modified version of the problem.

\begin{table}[!tbh]
\begin{center}
{\small
\begin{tabular}{lrlrrr}
\toprule
& & & & \multicolumn{1}{c}{\qIP} &  \multicolumn{1}{c}{\qBB} \\
                                                    & \multicolumn{1}{c}{Original} &\multicolumn{1}{c}{Fixed} & \multicolumn{1}{c}{M bound} & \multicolumn{1}{c}{Modified instance} & \multicolumn{1}{c}{detected}\\
\multicolumn{1}{c}{QP instance} & \multicolumn{1}{c}{Optimal Value} &\multicolumn{1}{c}{Variable}  &  \multicolumn{1}{c}{\eqref{eq:genQPbdM}} & \multicolumn{1}{c}{Optimal Value}
&  \multicolumn{1}{c}{unbounded dual}\\
\midrule
qp20\_10\_1\_1.mat	&	-13.189	&	$x_1$=0.4660	&	2.07E+04	&	-13.189	&	43-th  	\\
qp20\_10\_1\_2.mat	&	11.6662	&	$x_1$=1.0000	&	1.89E+04	&	11.6662	&	66-th  	\\
qp20\_10\_1\_4.mat	&	-18.3137	&	$x_1$=0.0000	&	1.29E+04	&	-18.3137	&	45-th  	\\
qp20\_10\_2\_1.mat	&	-3.2442	&	$x_1$=0.0000	&	6.98E+04	&	-3.2442	&	45-th  	\\
qp20\_10\_2\_2.mat	&	8.5919	&	$x_1$=0.0000	&	1.27E+04	&	8.5919	&	45-th  	\\
qp20\_10\_2\_4.mat	&	6.5794	&	$x_1$=0.0000	&	9.54E+03	&	6.5794	&	45-th  	\\
qp20\_10\_3\_1.mat	&	-30.179	&	$x_1$=0.0000	&	7.10E+04	&	-30.179	&	45-th  	\\
qp20\_10\_3\_2.mat	&	-15.0508	&	$x_1$=0.0000	&	4.70E+04	&	-15.0508	&	45-th  	\\
qp20\_10\_3\_4.mat	&	-12.665	&	$x_1$=0.0000	&	1.49E+04	&	-12.665	&	45-th  	\\
qp30\_15\_1\_1.mat	&	32.9577	&	$x_1$=0.0000	&	1.96E+05	&	32.9577	&	67-th  	\\
qp30\_15\_1\_3.mat	&	0.525	&	$x_1$=0.0000	&	3.91E+04	&	0.525	&	67-th  	\\
qp30\_15\_1\_4.mat	&	9.2296	&	$x_1$=0.0000	&	4.45E+04	&	9.2296	&	67-th  	\\
qp30\_15\_2\_3.mat	&	-2.0693	&	$x_1$=1.0000	&	3.18E+04	&	-2.0693	&	98-th  	\\
qp30\_15\_2\_4.mat	&	1.2862	&	$x_1$=0.0000	&	4.16E+04	&	1.2862	&	67-th  	\\
qp40\_20\_1\_3.mat	&	-2.7293	&	$x_1$=0.0000	&	8.30E+04	&	-2.7293	&	89-th  	\\
qp50\_25\_1\_4.mat	&	13.8442	&	$x_1$=0.0000	&	1.70E+05	&	13.8442	&	111-th  	\\
qp50\_25\_2\_4.mat	&	-6.8577	&	$x_1$=0.0000	&	2.80E+05	&	-6.8577	&	111-th  	\\
qp50\_25\_3\_2.mat	&	35.9871	&	$x_1$=0.0000	&	2.15E+05	&	35.9871	&	111-th  	\\
\bottomrule
\end{tabular}
}
\end{center}
\caption{Solution of QP test instances modified to have an unbounded dual feasible set using \qIP.}
\label{tab:unbounded}
\end{table}

\section{Computational results}
\label{sec:3}

In this section, we provide a detailed description of the implementation of the solution approach for  QP  problems described in the previous sections. Also, we illustrate the performance of the solution
approach by presenting the results of numerical experiments on a diverse set of  QP  test problems.

\subsection{Problem instances}
\label{sec:instances}

To test the performance of the proposed solution approach for  QP, we use the  set of
BoxQP \eqref{eq:BQP}, Globallib (\url{http://www.gamsworld.org/global/Globallib.htm}), CUTEr \citep{gould2003},
and RandQP test problems used in \cite[][Section 4.2 and Table 1]{chen2012}.
In addition to these test problems, we consider the following~ QP  test instances:

\begin{itemize}
\item SQP. Standard quadratic programming instances \eqref{eq:SQP} are created by replacing the constraints of each of the BoxQPs considered in~\cite[][Section 4.2 and Table 1]{chen2012}
by the constraint that the decision variables belong to the standard simplex of appropriate dimension. 

\item SQP30, SQP50 (see, \url{http://or.dei.unibo.it/library/msc}). A set of 300 SQP instances used for test purposes in~\cite{bonami2016solving}.
\item StableQP. These instances are particular SQPs resulting from the problem of computing the {\em stability number} of a graph~\citep[see, e.g.,][]{mot1965}. We use
instances of this type arising from a class of graphs that have been used for testing purposes in the literature  \citep[see, e.g.,][Section 4.2.2]{do2013}. A more detailed description of these instances is presented
in Section \ref{sec:graphs}.
\item Scozzari/Tardella (see, \url{http://or.dei.unibo.it/library/msc}). A set of 14 SQP instances used for test purposes in~\cite{scozzari2008clique, bonami2016solving}.
\item QPLIB2014 (see, \url{http://www.lamsade.dauphine.fr/QPlib2014/doku.php}). Nine nonconvex quadratic instances are selected from this test set. Four of the instances which are SQP instances are added to the SQP test set, and the other five instances are BoxQP instances, which are added to the BoxQP test set.
\end{itemize}

Similar to \citet{chen2012}, Table~\ref{table:instances} provides a summary of the basic information of all the test instances. In Table~\ref{table:instances},
$n$ denotes the range of the number of decision variables required to formulate the corresponding problem instance using $m_{ineq}$
inequality constraints, and $m_{eq}$ equality constraints. Also, {\it density} denotes the corresponding density range for the matrix defining the quadratic problem's objective.

\begin{table}[!htb]
\begin{center}
{\small
  \begin{tabular}{lrrrrr}
    \toprule
    Type & \# Instances & & \multicolumn{1}{c}{$n$} & \multicolumn{1}{c}{$m_{ineq}+m_{eq}$} & density \\
    \midrule
    StableQP	& 8	&& [5, 26] & [0,1] & [0.30, 0.60]\\
    SQP    & 90  && [20, 100] & [0, 90] & [0.19, 0.99]\\
    BoxQP    & 90  && [20, 100] & [0, 0] & [0.19, 0.99]\\
    Globallib & 83  && [2, 100] & [1,52] & [0.01, 1]\\
    CUTEr    & 6    && [4, 12] & [0, 13] & [0.08, 1]\\
    RandQP & 64  && [20, 50] & [14, 35] & [0.23, 1]\\
    SQP30 & 150 && [30, 30] & [0,1] & [1, 1]\\
    SQP50 & 150 && [50, 50] & [0,1] & [1, 1]\\
    Scozzari/Tardella & 14 && [30, 1000] & [0,1] & [0.25, 1]\\   
     \bottomrule
  \end{tabular}
  }
\end{center}
\caption{Statistics of the test  QP  instances.}
\label{table:instances}
\end{table}

\subsubsection{StableQP instances}
\label{sec:graphs}

For any graph $G$, the inverse of~$\alpha(G)$, the {\em stability number} of $G$, can be computed by solving the following SQP~\citep[see, e.g.,][]{mot1965}.
\begin{align}
\label{eq:stable}
\frac{1}{\alpha(G)} = \min_{x \in \Delta} x\tr (A+I)x,
\end{align}

where $A \in \S^n$ is the adjacency matrix of the undirected graph $G(V,E)$ with number of vertices
$\|V\| = n$, and set of edges $E \in V \times V$. Also $I$ is the identity matrix of appropriate
dimensions.

The StableQP instances are obtained by solving $\eqref{eq:stable}$ for a class
of graphs $G_k$, $k=1,\dots$ introduced in \cite{do2013}
that have proven to be hard instances for approximation methods for $\alpha(G)$
proposed in~\cite{bom2010,bun2009,dong2013, do2013}.

\subsection{Implementation details}
\label{sec:Implemetation}

The solution approach for  QP  proposed here is implemented as follows. First,
explicit upper and lower bounds for the instance's decision variables are obtained. Then,
the problem instance is reformulated as  QP  by linearly shifting its decision
variables, and adding slack variables to the problem as necessary (e.g., \eqref{eq:BQP_std}).
The upper bounds on the added slack variables are computed using \eqref{eq:primalBound}
to obtain the primal variable upper bounds~$U \in \R^n$.  Upper bounds $V \in \R^n$ on the dual variable are calculated using the methods
described in Section~\ref{sec:comphauff} (see \eqref{eq:sqpbd}, \eqref{eq:boxbd} and \eqref{eq:genbd}). Finally, CPLEX~12.5.1 (cf., \url{http://www-eio.upc.edu/lceio/manuals/CPLEX-11/html/})
is used to solve~IQP. The following parameter settings are used for CPLEX MILP solver:
\begin{itemize}
\item Max\_time: This is the user specified maximum running time of the algorithm and is set to $10^4$ seconds. Any problem taking longer than this value to be solved will be deemed as ``out of time''.
\item Tol: The solver will stop when
\[
\frac{|\text{bestnode}-\text{bestinteger}|}{1^{-10}+|\text{bestinteger}| } \le 10^{-6}.
\]
For the interested reader, the definition of the parameters {\em bestnode} and {\em bestinteger} can be found in~\cite{cplex201012}. Here, it suffices to say that this criteria is consistent with {\tt quadprogBB} stopping criteria~\citep[cf.,][]{chen2012}, which is
\[
\frac{\text{Greatest upper bound} - \text{current lower bound}}{max\{1,|\text{Greatest upper bound}| \} }\le 10^{-6}.
\]
\item Other parameters of the CPLEX MILP solver such as TolXInteger, Max\_iter, BranchStrategy, Nodeselect, are set to their default values.
\end{itemize}

We refer to the procedure described in this section to solve  QP  as {\tt quadprogIP}, which is coded using
{\tt Matlab R2014a}, and is publicly available at \url{https://github.com/xiawei918/quadprogIP}.

\subsection{Numerical performance}

In order to test the performance of the~{\tt quadprogIP} methodology proposed in Section~\ref{sec:Implemetation}, the QP test instances discussed in
Section~\ref{sec:instances} are solved using {\tt quadprogIP}, the {\tt quadprogBB} solver introduced by~\citet{chen2012}, the NLP solver \BARON{} {\tt 17.8.9} of ~\citet{Sahi96}, and the {\tt CPLEX 12.7.0.0} QP solver.
All tests are done using {\tt Matlab R2014b (8.4.0.150421)}, together with CPLEX~12.7.0.0.,  on a
AMD Opteron 2.0 GHz  machine with 32GB memory and
16 cores (each core is a 2.0 GHz. 64 bit architecture),
from the COR@L laboratory (cf., \url{http://coral.ise.lehigh.edu/}).

Similar to \citet{chen2012}, to compare the performance between {\tt quadprogIP} and {\tt quadprogBB}, \qIP and \BARON, and \qIP and \cplex, we plot the solution time it takes to solve a particular  QP  instance with two of the solvers as a square
in a 2D plane, where the $y$-axis denotes either {\tt quadprogBB}'s, \BARON{}'s, or \cplex's  solution time and the $x$-axis
 denotes {\tt quadprogIP} 's solution time.
The dashed line in the plots indicates the $y=x$ line in the plane, that represents equal colution solution times. Thus, a square that is above the diagonal line indicates an instance for which it takes the solver represented on $y$-axis more solution time to solve than {\tt quadprogIP}.
Furthermore, the size of the square illustrates the size (number of decision variables) of the instance. That is, smaller squares represent ``smaller'' size
instances while bigger squares represent ``bigger'' size instances.
In the figures below, only instances in which at least one of the methodologies solves
the problem within the maximum time allowed are displayed.

\subsubsection{Results on SQP instances.}
\label{sec:SQPnum}
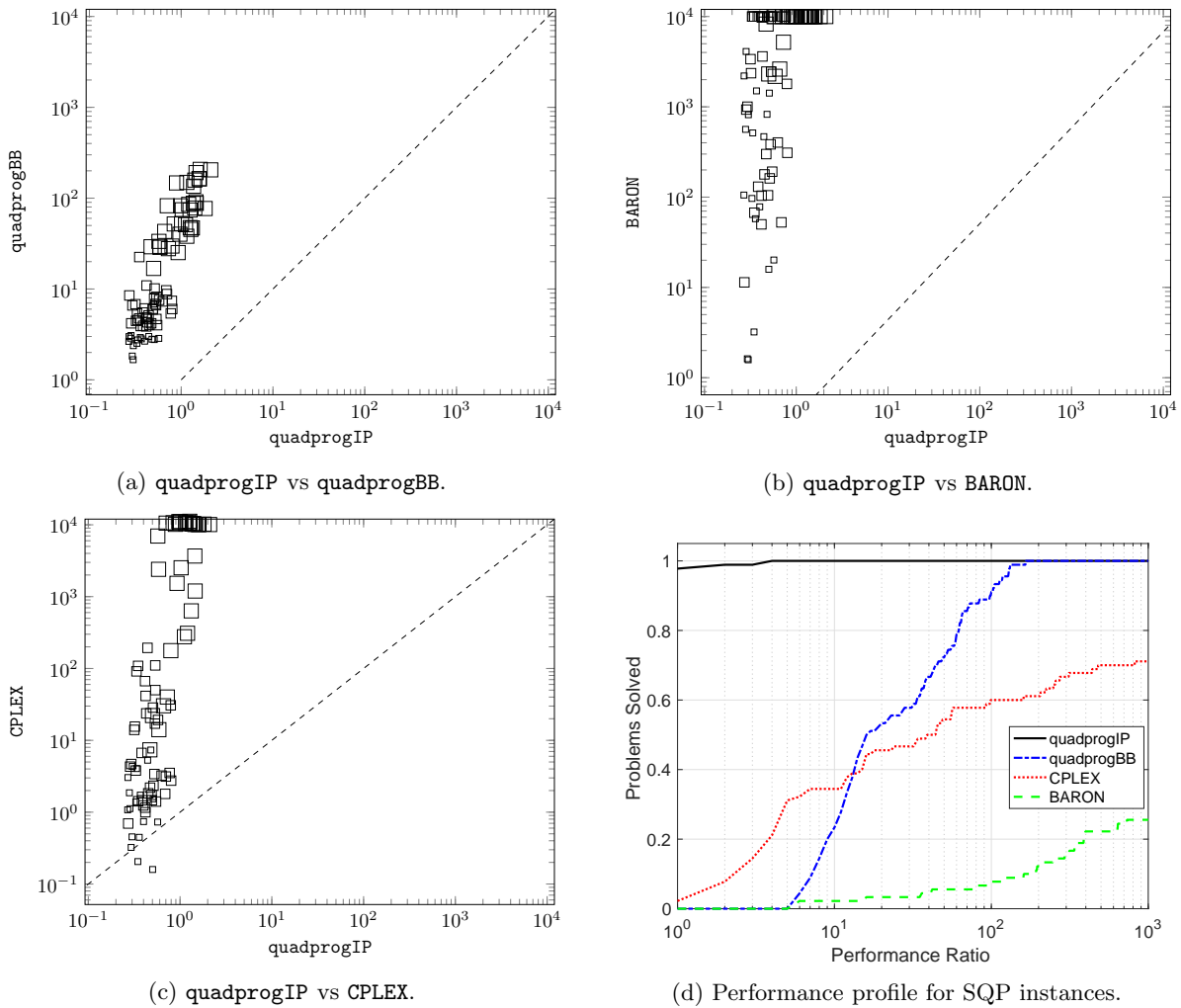
\begin{figure}[!htb]
    \captionsetup[figure]{justification=left}
    \captionsetup[subfigure]{justification=centering}
        \centering
        \begin{subfigure}{0.495\textwidth}
            \centering
            \begin{tikzpicture}[scale = \myplotscale]
	\begin{loglogaxis}[%
	height = 8.5cm,
	width = 10cm,
	xmin=0,
  xmax=12000,
  ymin=0,
  ymax=12000,
    xlabel={{\tt quadprogIP}},
    ylabel={{\tt quadprogBB}},
	scatter/classes={%
		a={mark=square,black, scale = 0.75},
		b={mark=square,black, scale = 1.2},
		c={mark=square,black, scale = 1.8}}]
	\addplot[scatter,only marks,%
		scatter src=explicit symbolic]%
	table[meta=label] {
x     y      label
0.3468	3.7408	a
0.2996	1.6665	a
0.2931	1.8271	a
0.3293	2.5231	a
0.4002	2.6651	a
0.3010	2.3780	a
0.2689	2.6502	a
0.3359	2.7509	a
0.4806	2.8006	a
0.2805	2.8729	a
0.4437	2.9994	a
0.2702	2.9921	a
0.2843	3.0762	a
0.3688	2.8232	a
0.5091	2.7736	a
0.3590	2.9220	a
0.5027	6.8565	a
0.5703	2.8585	a
0.2723	8.4928	b
0.4194	10.9331	b
0.3848	3.9363	b
0.4173	3.8720	b
0.5474	3.9978	b
0.4503	4.0162	b
0.4713	4.1920	b
0.2855	4.1915	b
0.2941	6.5934	b
0.3220	4.4646	b
0.3174	6.8007	b
0.4156	4.6633	b
0.4285	4.2967	b
0.4217	4.8493	b
0.3492	4.7473	b
0.5370	4.6774	b
0.3352	4.5432	b
0.4418	5.1192	b
0.3810	5.5115	b
0.3301	5.3821	b
0.7702	5.3897	b
0.7005	8.8003	b
0.3983	6.1081	b
0.4901	5.9507	b
0.4925	6.0527	b
0.5257	6.8548	b
0.7997	6.0228	b
0.6296	7.2951	b
0.5754	7.5829	b
0.5326	6.6261	b
0.7919	7.3981	b
0.5072	7.8206	b
0.5381	8.1969	b
0.6889	9.6652	b
0.5137	10.1314	b
0.3497	22.4811	b
0.5879	28.7480	c
0.4693	29.1757	c
0.5011	16.8481	c
0.5833	29.4741	c
0.9264	25.1653	c
0.5724	33.4560	c
1.2557	47.4539	c
0.8399	52.0111	c
1.1516	38.1031	c
0.7941	29.7090	c
0.6611	43.0305	c
0.7263	27.8134	c
0.9688	40.1076	c
1.2759	45.3056	c
1.0094	52.5040	c
1.3534	87.9029	c
1.4052	77.9549	c
1.8243	76.9567	c
1.1155	51.0390	c
1.3245	47.0459	c
1.2066	85.4548	c
0.7053	82.5460	c
1.2486	74.1158	c
0.9859	82.4277	c
1.5637	164.0188	c
0.8894	147.0859	c
1.5853	162.8108	c
1.4608	87.6446	c
1.4484	91.1178	c
1.0267	74.5291	c
1.3661	135.4834	c
1.4060	157.9956	c
1.1530	149.1907	c
1.6110	207.6056	c
1.4531	191.1164	c
2.1002	204.9668	c
};
\draw [dashed] (0,0) -- (120,120);
	\end{loglogaxis}
\end{tikzpicture}
            \caption{\qIP\ vs \qBB.\vspace{0.3\baselineskip} \quad}    
            \label{fig:SQP_BB}
        \end{subfigure}
        \begin{subfigure}{0.495\textwidth}  
            \centering 
            \begin{tikzpicture}[scale = \myplotscale]
	\begin{loglogaxis}[%
	height = 8.5cm,
	width = 10cm,
	xmin=0,
  xmax=12000,
  ymin=0,
  ymax=12000,
    xlabel={{\tt quadprogIP}},
    ylabel={{\tt \BARON{}}},
	scatter/classes={%
		a={mark=square,black, scale = 0.75},
		b={mark=square,black, scale = 1.2},
		c={mark=square,black, scale = 1.8}}]
	\addplot[scatter,only marks,%
		scatter src=explicit symbolic]%
	table[meta=label] {
x     y      label
0.3468	3.2100	a
0.2996	1.5800	a
0.2931	1.6200	a
0.3293	96.7100	a
0.4002	78.0800	a
0.3010	817.6500	a
0.2689	105.4300	a
0.3359	514.1100	a
0.4806	826.1300	a
0.2805	562.9100	a
0.4437	467.3600	a
0.2702	2199.8100	a
0.2843	4102.6500	a
0.3688	1501.6900	a
0.5091	1419.1200	a
0.3590	57.6400	a
0.5027	15.8700	a
0.5703	20.1200	a
0.2723	11.3800	b
0.4194	103.7300	b
0.3848	130.4400	b
0.4173	49.9300	b
0.5474	191.3700	b
0.4503	178.6000	b
0.4713	299.8000	b
0.2855	932.0700	b
0.2941	1003.6200	b
0.3220	2359.3900	b
0.3174	3372.8700	b
0.4156	10000.0000	b
0.4285	3615.4500	b
0.4217	9762.2000	b
0.3492	10000.0400	b
0.5370	10000.0700	b
0.3352	10000.0500	b
0.4418	10000.0300	b
0.3810	10000.0400	b
0.3301	10000.0700	b
0.7702	10000.1100	b
0.7005	10000.0200	b
0.3983	10000.0300	b
0.4901	10000.0700	b
0.4925	104.7700	b
0.5257	384.2200	b
0.7997	309.8300	b
0.6296	402.6000	b
0.5754	10000.0500	b
0.5326	2448.0800	b
0.7919	1796.1200	b
0.5072	10000.0200	b
0.5381	10000.0800	b
0.6889	52.5400	b
0.5137	161.2600	b
0.3497	67.2100	b
0.5879	2177.7200	c
0.4693	8256.4500	c
0.5011	2313.2800	c
0.5833	10000.0500	c
0.9264	10000.7600	c
0.5724	10000.0600	c
1.2557	10000.0100	c
0.8399	10000.1800	c
1.1516	10000.0000	c
0.7941	10000.3400	c
0.6611	2620.3800	c
0.7263	5193.8300	c
0.9688	10000.0000	c
1.2759	10000.0500	c
1.0094	10000.0400	c
1.3534	10000.0100	c
1.4052	10000.5100	c
1.8243	10000.0800	c
1.1155	10000.1300	c
1.3245	10000.0600	c
1.2066	10000.0500	c
0.7053	10000.0100	c
1.2486	10000.1000	c
0.9859	10000.0200	c
1.5637	10000.0700	c
0.8894	10000.2600	c
1.5853	10000.1800	c
1.4608	10000.3200	c
1.4484	10000.2700	c
1.0267	10000.0500	c
1.3661	10000.3400	c
1.4060	10000.3200	c
1.1530	10000.2100	c
1.6110	10000.3300	c
1.4531	10000.6200	c
2.1002	10000.0400	c
};
\draw [dashed] (-8.5,-10) -- (95,100);
	\end{loglogaxis}
\end{tikzpicture}
            \caption[]%
            {\qIP\ vs \BARON.\vspace{0.3\baselineskip}  \quad}    
            \label{fig:SQP_BA}
        \end{subfigure}
        \begin{subfigure}{0.495\textwidth}   
            \centering 
            \begin{tikzpicture}[scale = \myplotscale]
	\begin{loglogaxis}[%
	height = 8.5cm,
	width = 10cm,
	xmin=0,
  xmax=12000,
  ymin=0,
  ymax=12000,
    xlabel={{\tt quadprogIP}},
    ylabel={{\tt CPLEX}},
	scatter/classes={%
		a={mark=square,black, scale = 0.75},
		b={mark=square,black, scale = 1.2},
		c={mark=square,black, scale = 1.8}}]
	\addplot[scatter,only marks,%
		scatter src=explicit symbolic]%
	table[meta=label] {
x     y      label
0.3468	0.2054	a
0.2996	0.4537	a
0.2931	0.3232	a
0.3293	1.3905	a
0.4002	0.7401	a
0.3010	4.3389	a
0.2689	1.0798	a
0.3359	4.0476	a
0.4806	7.3341	a
0.2805	1.8624	a
0.4437	5.2863	a
0.2702	3.0368	a
0.2843	1.1199	a
0.3688	1.7276	a
0.5091	1.5358	a
0.3590	0.4470	a
0.5027	0.1594	a
0.5703	0.7290	a
0.2723	0.6955	b
0.4194	0.9983	b
0.3848	1.4651	b
0.4173	1.1570	b
0.5474	1.4051	b
0.4503	1.8312	b
0.4713	2.1933	b
0.2855	4.3053	b
0.2941	4.6808	b
0.3220	15.5413	b
0.3174	13.9872	b
0.4156	66.6155	b
0.4285	23.7419	b
0.4217	41.5044	b
0.3492	108.6921	b
0.5370	110.3878	b
0.3352	91.2234	b
0.4418	193.6757	b
0.3810	6.6831	b
0.3301	3.7998	b
0.7702	3.5151	b
0.7005	3.1339	b
0.3983	1.2421	b
0.4901	1.7427	b
0.4925	1.4676	b
0.5257	3.3965	b
0.7997	2.7502	b
0.6296	3.2307	b
0.5754	19.2290	b
0.5326	17.0766	b
0.7919	30.7508	b
0.5072	28.5458	b
0.5381	50.0078	b
0.6889	1.7990	b
0.5137	2.3344	b
0.3497	1.4182	b
0.5879	14.0787	c
0.4693	7.3794	c
0.5011	22.1389	c
0.5833	2399.6981	c
0.9264	1531.7942	c
0.5724	6975.7964	c
1.2557	10520.8018	c
0.8399	10682.6410	c
1.1516	10705.2334	c
0.7941	177.1826	c
0.6611	30.2729	c
0.7263	40.0808	c
0.9688	11073.9027	c
1.2759	11199.8436	c
1.0094	11062.7860	c
1.3534	10272.7751	c
1.4052	10286.4725	c
1.8243	10227.4466	c
1.1155	276.7451	c
1.3245	633.5925	c
1.2066	309.1827	c
0.7053	10548.8333	c
1.2486	10731.3723	c
0.9859	10624.9357	c
1.5637	10140.4094	c
0.8894	10155.5781	c
1.5853	10138.6014	c
1.4608	1197.1191	c
1.4484	3676.3774	c
1.0267	2522.3906	c
1.3661	10280.9034	c
1.4060	10358.7360	c
1.1530	10334.1443	c
1.6110	10096.2472	c
1.4531	10092.9874	c
2.1002	10081.6796	c
};
\draw [dashed] (-10,-10) -- (120,120);
	\end{loglogaxis}
\end{tikzpicture}
            \caption[]%
            {\qIP\ vs \cplex.\vspace{0.3\baselineskip}  \quad}    
            \label{fig:SQP_CP}
        \end{subfigure}
        \begin{subfigure}{0.495\textwidth}   
            \centering 
            \includegraphics[width=\textwidth]{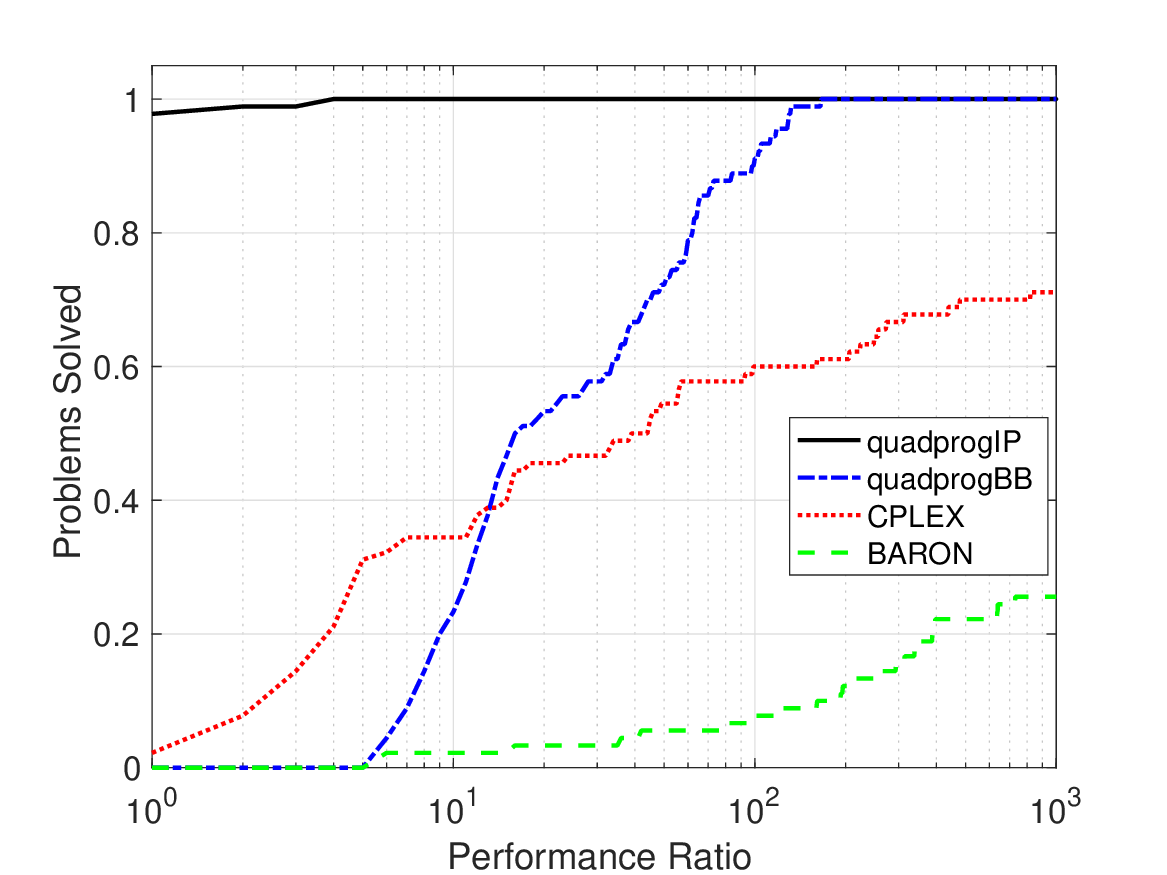}
            \caption[]%
            {Performance profile for SQP instances.\vspace{0.3\baselineskip}  \quad}    
            \label{fig:SQP_Profile}
        \end{subfigure}
        \caption[]
        {Solution time in seconds of SQP instances. Size of squares illustrates size of the instance. A square at the maximum value of an axis represents an instances for which the solver in that axis reached maximum running time without solving it.} 
        \label{fig:SQP_ALL}
    \end{figure}

The results for the SQP test instances are shown in Figure~\ref{fig:SQP_ALL}. Note that a different scale is used in the axes of Figures~\ref{fig:SQP_BB},~\ref{fig:SQP_BA}, and~\ref{fig:SQP_CP}.

 Figure~\ref{fig:SQP_BB} shows that \qIP  clearly outperforms \qBB\ by solving all SQP instances in a time that is one to two orders of magnitude faster than \qBB\, and specially in the larger instances. Similarly, Figure~\ref{fig:SQP_BA} shows that \qIP  clearly outperforms \BARON by solving all SQP instances in a time that is one to two orders of magnitude faster than \BARON, and specially in the larger instances. Although \cplex solves two small-scale instances faster than \qIP, again, in general \qIP outperforms \cplex by orders of magnitude in terms of solution time  (see, Figure~\ref{fig:SQP_CP}).  
As Figures~\ref{fig:SQP_BB},~\ref{fig:SQP_BA}, and~\ref{fig:SQP_CP} illustrate, the performance of \qIP\  against the other solvers improves as the SQP instance becomes larger.
The performance profile in Figure~\ref{fig:SQP_Profile} summarizes the clear advantages of solving the very important class of SQP instances with the proposed \qIP solution approach.

\subsubsection{Results on SQP30 and SQP50 instances.}
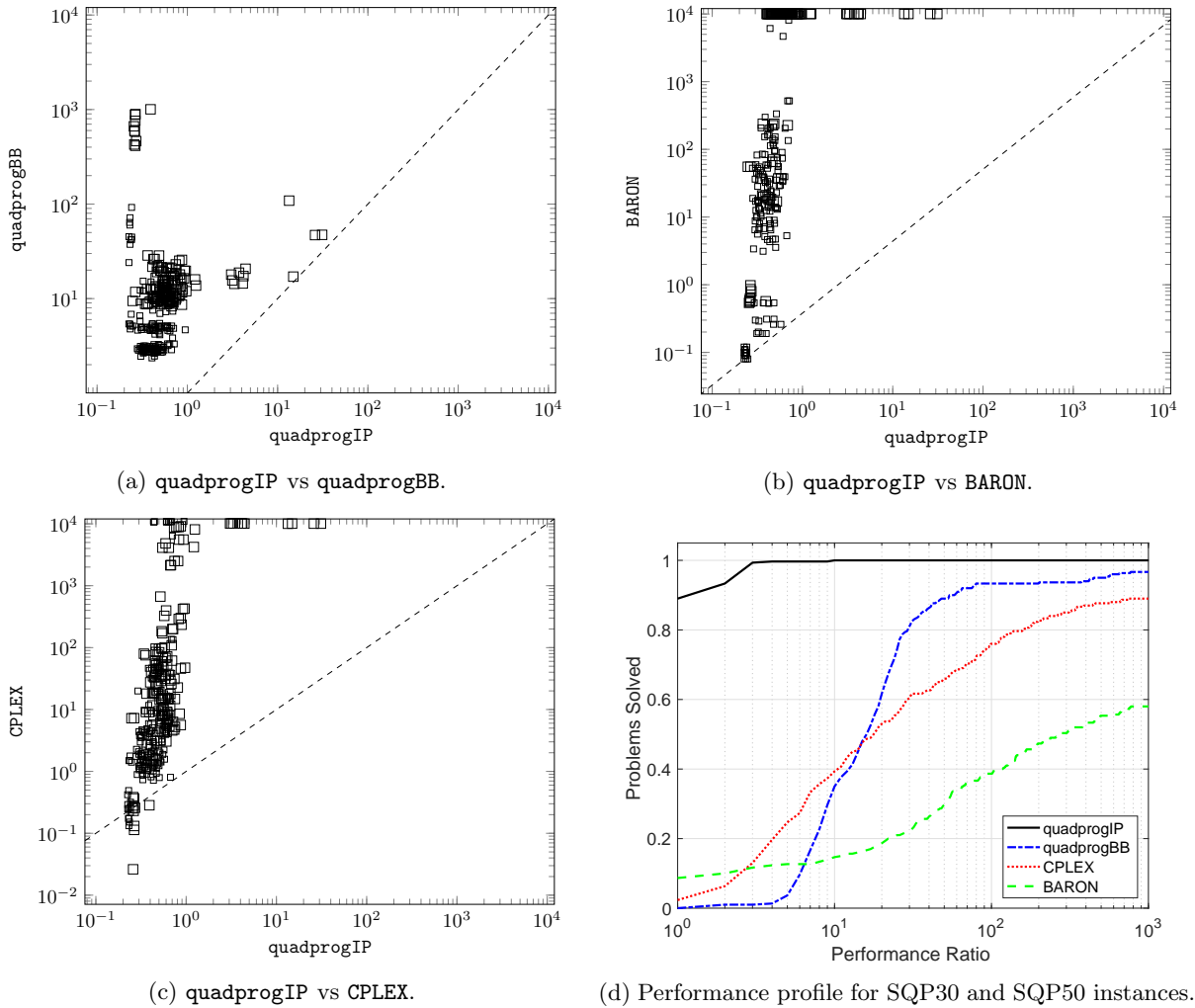
\begin{figure}[!htb]
    \captionsetup[figure]{justification=left}
    \captionsetup[subfigure]{justification=centering}
        \centering
        \begin{subfigure}{0.495\textwidth}
            \centering
            \begin{tikzpicture}[scale = \myplotscale]
	\begin{loglogaxis}[%
	height = 8.5cm,
	width = 10cm,
	xmin=0,
  xmax=12000,
  ymin=0,
  ymax=12000,
    xlabel={{\tt quadprogIP}},
    ylabel={{\tt quadprogBB}},
	scatter/classes={%
		a={mark=square,black, scale = 0.75},
		b={mark=square,black, scale = 1.2},
		c={mark=square,black, scale = 1.8}}]
	\addplot[scatter,only marks,%
		scatter src=explicit symbolic]%
	table[meta=label] {
x     y      label
0.4975	6.9541	a
0.3817	4.7952	a
0.4375	4.2000	a
0.4396	4.5897	a
0.4778	4.5490	a
0.6686	2.8744	a
0.6939	8.7935	a
0.4941	5.4671	a
0.6090	4.6021	a
0.4263	3.2608	a
0.4109	5.3668	a
0.6045	3.3228	a
0.6384	3.3099	a
0.6116	3.1842	a
0.6055	4.8926	a
0.6860	3.8308	a
0.7164	3.3145	a
0.7044	4.8467	a
0.4895	4.2943	a
0.4960	4.1432	a
0.6905	7.3608	a
0.5285	8.2963	a
0.4506	4.9447	a
0.4459	2.9414	a
0.4381	5.1108	a
0.4150	4.9064	a
0.9468	4.6712	a
0.4784	3.2752	a
0.4409	2.9008	a
0.6390	4.8241	a
0.4146	7.5432	a
0.5185	3.0681	a
0.4691	2.7758	a
0.3999	2.9129	a
0.4415	3.0134	a
0.5341	2.7333	a
0.3426	4.9779	a
0.4301	2.9134	a
0.3838	3.1072	a
0.5270	2.8918	a
0.4854	5.0826	a
0.4447	2.6943	a
0.4456	3.0885	a
0.4147	3.1244	a
0.4710	2.7940	a
0.4316	2.8592	a
0.3194	3.0137	a
0.3473	2.9166	a
0.2818	4.6319	a
0.3804	2.8995	a
0.6049	5.1328	a
0.5007	3.1959	a
0.4094	5.4382	a
0.5149	5.1323	a
0.3279	3.1562	a
0.3022	3.1268	a
0.3861	2.9026	a
0.4358	2.8335	a
0.4189	2.7885	a
0.3047	2.4527	a
0.4909	5.0050	a
0.3426	3.1579	a
0.5152	2.9170	a
0.3784	2.9402	a
0.2976	6.5697	a
0.5979	5.1104	a
0.3595	2.8853	a
0.4454	2.8990	a
0.4638	2.8324	a
0.4499	2.8490	a
0.5686	5.0884	a
0.6415	4.5311	a
0.3709	2.9835	a
0.3317	3.0427	a
0.3193	2.9071	a
0.3890	4.5871	a
0.3297	3.0105	a
0.5410	3.0495	a
0.2989	2.9431	a
0.3864	3.0814	a
0.5938	4.9608	a
0.3289	2.8793	a
0.3130	2.9462	a
0.3943	2.9435	a
0.3256	2.6522	a
0.3734	2.7729	a
0.4095	2.9922	a
0.5136	2.6965	a
0.3483	2.7072	a
0.4060	2.3582	a
0.3052	4.9181	a
0.2994	5.7915	a
0.3193	2.8194	a
0.3336	2.7838	a
0.2869	2.8444	a
0.3745	2.7510	a
0.3105	4.7009	a
0.6704	5.1271	a
0.3138	2.8132	a
0.3090	2.7957	a
0.5048	4.8668	a
0.3341	2.8381	a
0.4274	2.6406	a
0.3268	2.7213	a
0.4323	2.4663	a
0.5004	2.8495	a
0.4342	2.8047	a
0.4084	2.7538	a
0.4338	2.6845	a
0.4088	4.7070	a
0.3662	4.7212	a
0.3180	3.0283	a
0.3826	2.9632	a
0.5515	2.8589	a
0.3351	2.7405	a
0.3175	4.9203	a
0.3629	2.7039	a
0.4690	2.9260	a
0.4677	3.0081	a
0.2791	2.9472	a
0.2388	6.8304	a
0.2414	42.2799	a
0.2325	4.7744	a
0.4004	5.1246	a
0.2256	24.0475	a
0.4837	5.0631	a
0.2365	5.1282	a
0.2310	59.9543	a
0.2483	5.1294	a
0.3962	15.0211	a
0.2413	92.1966	a
0.3240	15.3150	a
0.3677	9.7031	a
0.2294	63.9459	a
0.4230	9.7135	a
0.4825	9.0810	a
0.2420	44.7685	a
0.2919	9.2184	a
0.4771	5.0589	a
0.2319	37.4881	a
0.5767	7.2617	a
0.3005	12.2521	a
0.2357	41.8289	a
0.3288	12.0542	a
0.2309	5.3767	a
0.2320	71.7053	a
0.2247	5.3924	a
0.3118	4.9606	a
0.2253	45.2104	a
0.3588	4.9674	a
0.9329	19.6881	b
0.8337	16.6846	b
4.1844	17.3052	b
0.6683	18.6532	b
0.6790	20.7680	b
4.1127	14.4512	b
0.9613	19.6419	b
0.9209	16.3840	b
3.0645	17.9459	b
0.8144	25.1711	b
0.7441	23.0059	b
3.3178	14.3303	b
0.8256	17.9714	b
0.7794	15.2988	b
4.4178	20.6909	b
0.6432	21.2609	b
0.5462	19.1001	b
3.7384	18.8071	b
0.7037	17.5986	b
0.7170	15.6488	b
3.1523	15.5590	b
1.2164	15.9591	b
1.2485	13.7330	b
31.1552	47.2579	b
0.8713	25.6869	b
13.4146	108.7284	b
14.9197	17.0845	b
0.5918	20.7958	b
0.7356	18.5416	b
25.9848	46.9670	b
0.5340	14.2138	b
0.4001	12.6307	b
0.8275	11.1133	b
0.5547	11.7072	b
0.4433	13.5128	b
0.5355	15.2820	b
0.6157	11.6668	b
0.5353	12.9852	b
0.5723	14.8605	b
0.4247	12.3657	b
0.8469	13.3401	b
0.4950	11.8155	b
0.5373	15.6851	b
0.4464	16.5015	b
0.5460	10.8823	b
0.5244	13.4708	b
0.4621	19.1011	b
0.8702	15.9862	b
0.9726	11.9860	b
0.6979	17.3170	b
0.5182	19.7709	b
0.5292	12.2282	b
0.4901	12.4468	b
0.6915	11.3879	b
0.4844	22.5535	b
0.6172	21.1922	b
0.8885	12.4852	b
0.8747	11.6154	b
0.7515	12.7535	b
0.6005	15.5060	b
0.5450	13.0191	b
0.6686	14.7076	b
0.5442	10.4026	b
0.6333	11.1732	b
0.5831	11.5911	b
0.6157	14.2065	b
0.4218	10.3907	b
0.5242	10.7365	b
0.4514	9.9338	b
0.5321	10.3000	b
0.5512	12.3465	b
0.5921	9.7652	b
0.5857	10.4670	b
0.5615	11.3897	b
0.5537	10.3739	b
0.4100	21.7422	b
0.6200	10.7696	b
0.4574	10.7104	b
0.5824	10.3059	b
0.5834	11.5432	b
0.4646	10.4862	b
0.5583	13.9801	b
0.7534	10.1572	b
0.5359	10.1963	b
0.8077	10.6782	b
0.5470	11.0257	b
0.5365	10.9408	b
0.6513	11.5065	b
0.6230	15.4454	b
0.4730	10.3068	b
0.5084	8.7892	b
0.6013	12.2368	b
0.6687	8.2995	b
0.7157	10.4405	b
0.5481	15.9039	b
0.5372	13.2282	b
0.4341	9.8401	b
0.7381	14.0534	b
0.7409	9.7482	b
0.6652	9.9617	b
0.6207	10.7819	b
0.4577	9.2108	b
0.7014	9.0816	b
0.6515	9.3084	b
0.6397	9.4543	b
0.5356	9.6352	b
0.5798	10.7076	b
0.4970	10.1488	b
0.8712	8.6167	b
0.6035	9.7468	b
0.5963	8.7796	b
0.6082	9.5305	b
0.4083	10.5890	b
0.6196	8.9864	b
0.6309	9.5616	b
0.5674	13.3031	b
0.6298	9.1616	b
0.6366	9.4667	b
0.5177	9.5776	b
0.5584	12.1846	b
0.2650	11.7463	b
0.2715	463.7084	b
0.2479	9.4518	b
0.4861	28.5215	b
0.2620	587.1670	b
0.3636	28.5930	b
0.4337	26.1701	b
0.2585	427.1688	b
0.4213	26.2187	b
0.3581	12.6286	b
0.2632	881.5047	b
0.3900	12.7512	b
0.3555	8.7236	b
0.2652	413.9255	b
0.3407	8.8045	b
0.4621	21.9931	b
0.2653	759.5818	b
0.5027	19.7773	b
0.3680	8.9328	b
0.3916	1005.9033	b
0.4101	8.8404	b
0.5199	21.4501	b
0.2593	595.6642	b
0.4770	21.4149	b
0.5783	13.0504	b
0.2663	894.1126	b
0.5391	13.1467	b
0.6866	9.3720	b
0.2558	664.7486	b
0.4723	9.3784	b
};
\draw [dashed] (0,0) -- (120,120);
	\end{loglogaxis}
\end{tikzpicture}
            \caption{\qIP\ vs \qBB.\vspace{0.3\baselineskip} \quad}    
            \label{fig:SQP35_BB}
        \end{subfigure}
        \begin{subfigure}{0.495\textwidth}  
            \centering 
            \begin{tikzpicture}[scale = \myplotscale]
	\begin{loglogaxis}[%
	height = 8.5cm,
	width = 10cm,
	xmin=0,
  xmax=12000,
  ymin=0,
  ymax=12000,
    xlabel={{\tt quadprogIP}},
    ylabel={{\tt \BARON{}}},
	scatter/classes={%
		a={mark=square,black, scale = 0.75},
		b={mark=square,black, scale = 1.2},
		c={mark=square,black, scale = 1.8}}]
	\addplot[scatter,only marks,%
		scatter src=explicit symbolic]%
	table[meta=label] {
x     y      label
0.4975	153.2600	a
0.3817	153.3200	a
0.4375	6108.4100	a
0.4396	213.9300	a
0.4778	213.7400	a
0.6686	10000.1100	a
0.6939	134.7900	a
0.4941	134.9000	a
0.6090	4669.2800	a
0.4263	162.3900	a
0.4109	162.3600	a
0.6045	10000.0000	a
0.6384	39.3700	a
0.6116	39.3500	a
0.6055	10000.0900	a
0.6860	519.5300	a
0.7164	519.4800	a
0.7044	8103.2700	a
0.4895	95.9600	a
0.4960	95.9400	a
0.6905	10000.0800	a
0.5285	75.0600	a
0.4506	75.1100	a
0.4459	10000.0000	a
0.4381	204.8300	a
0.4150	204.9100	a
0.9468	10000.0200	a
0.4784	113.3600	a
0.4409	113.3900	a
0.6390	10000.0100	a
0.4146	58.5200	a
0.5185	13.5700	a
0.4691	142.2800	a
0.3999	12.4600	a
0.4415	7.9000	a
0.5341	48.0200	a
0.3426	48.9300	a
0.4301	9.0000	a
0.3838	83.4200	a
0.5270	13.0300	a
0.4854	8.9700	a
0.4447	13.8200	a
0.4456	24.4200	a
0.4147	13.9600	a
0.4710	123.6300	a
0.4316	36.4000	a
0.3194	52.0500	a
0.3473	207.6100	a
0.2818	18.8800	a
0.3804	18.4000	a
0.6049	33.4800	a
0.5007	14.0300	a
0.4094	28.0600	a
0.5149	334.4900	a
0.3279	16.1700	a
0.3022	17.0000	a
0.3861	298.0000	a
0.4358	12.1800	a
0.4189	14.8100	a
0.3047	59.3400	a
0.4909	8.0800	a
0.3426	11.4000	a
0.5152	35.1300	a
0.3784	34.7900	a
0.2976	6.4700	a
0.5979	73.8900	a
0.3595	85.0600	a
0.4454	6.4600	a
0.4638	61.7400	a
0.4499	23.7700	a
0.5686	29.5300	a
0.6415	204.7800	a
0.3709	8.7600	a
0.3317	6.9600	a
0.3193	17.2900	a
0.3890	32.3800	a
0.3297	6.5600	a
0.5410	55.5700	a
0.2989	36.2700	a
0.3864	18.0700	a
0.5938	90.2400	a
0.3289	42.0100	a
0.3130	10.1500	a
0.3943	58.8700	a
0.3256	16.6000	a
0.3734	6.7400	a
0.4095	52.5400	a
0.5136	22.9200	a
0.3483	10.2200	a
0.4060	22.3900	a
0.3052	12.5600	a
0.2994	8.7500	a
0.3193	28.3400	a
0.3336	55.1000	a
0.2869	3.3800	a
0.3745	37.8600	a
0.3105	44.8100	a
0.6704	5.3000	a
0.3138	82.2900	a
0.3090	9.9100	a
0.5048	3.5600	a
0.3341	8.0400	a
0.4274	16.8700	a
0.3268	13.3000	a
0.4323	31.0000	a
0.5004	4.7200	a
0.4342	4.3200	a
0.4084	33.5500	a
0.4338	6.8900	a
0.4088	4.6500	a
0.3662	3.1200	a
0.3180	17.6600	a
0.3826	5.1400	a
0.5515	28.9400	a
0.3351	13.1500	a
0.3175	5.5900	a
0.3629	42.9400	a
0.4690	19.5900	a
0.4677	4.7100	a
0.2791	52.6500	a
0.2388	0.0900	a
0.2414	0.1000	a
0.2325	0.0900	a
0.4004	0.3100	a
0.2256	0.1200	a
0.4837	0.3100	a
0.2365	0.0800	a
0.2310	0.1100	a
0.2483	0.0800	a
0.3962	0.1900	a
0.2413	0.1200	a
0.3240	0.2000	a
0.3677	0.5400	a
0.2294	0.1100	a
0.4230	0.5300	a
0.4825	0.5400	a
0.2420	0.1000	a
0.2919	0.5400	a
0.4771	0.2600	a
0.2319	0.1000	a
0.5767	0.2600	a
0.3005	0.3000	a
0.2357	0.1000	a
0.3288	0.2900	a
0.2309	0.0900	a
0.2320	0.1100	a
0.2247	0.1000	a
0.3118	0.1900	a
0.2253	0.0800	a
0.3588	0.1900	a
		
0.9329	10000.0800	b
0.8337	10000.0600	b
4.1844	10000.0700	b
0.6683	10000.0600	b
0.6790	10000.0000	b
4.1127	10000.1700	b
0.9613	10000.1500	b
0.9209	10000.0000	b
3.0645	10000.1100	b
0.8144	10000.0500	b
0.7441	10000.0500	b
3.3178	10000.1600	b
0.8256	10000.0400	b
0.7794	10000.0700	b
4.4178	10000.0200	b
0.6432	10000.0100	b
0.5462	10000.0100	b
3.7384	10000.1600	b
0.7037	10000.0400	b
0.7170	10000.0400	b
3.1523	10000.0000	b
1.2164	10000.0800	b
1.2485	10000.0200	b
31.1552	10000.1600	b
0.8713	10000.0200	b
13.4146	10000.0900	b
14.9197	10000.1600	b
0.5918	10000.1500	b
0.7356	10000.0400	b
25.9848	10000.0100	b
0.5340	10000.0100	b
0.4001	10000.0900	b
0.8275	10000.0900	b
0.5547	10000.1000	b
0.4433	10000.0900	b
0.5355	10000.1200	b
0.6157	10000.1000	b
0.5353	10000.0200	b
0.5723	10000.0300	b
0.4247	10000.1000	b
0.8469	10000.0600	b
0.4950	10000.0400	b
0.5373	10000.0000	b
0.4464	10000.1000	b
0.5460	10000.0700	b
0.5244	10000.1100	b
0.4621	10000.0500	b
0.8702	10000.1200	b
0.9726	10000.0400	b
0.6979	10000.0400	b
0.5182	10000.0000	b
0.5292	10000.0600	b
0.4901	10000.1000	b
0.6915	10000.0100	b
0.4844	10000.1100	b
0.6172	10000.0000	b
0.8885	10000.1000	b
0.8747	10000.1100	b
0.7515	10000.0000	b
0.6005	10000.1100	b
0.5450	10000.1600	b
0.6686	10000.0700	b
0.5442	10000.0500	b
0.6333	10000.1000	b
0.5831	10000.0700	b
0.6157	10000.1000	b
0.4218	10000.0000	b
0.5242	10000.0100	b
0.4514	10000.0400	b
0.5321	10000.0900	b
0.5512	10000.0200	b
0.5921	10000.0200	b
0.5857	10000.0200	b
0.5615	10000.0900	b
0.5537	10000.0000	b
0.4100	10000.0700	b
0.6200	10000.2000	b
0.4574	10000.0700	b
0.5824	10000.0000	b
0.5834	10000.0400	b
0.4646	10000.0900	b
0.5583	10000.0900	b
0.7534	10000.0900	b
0.5359	10000.0500	b
0.8077	10000.0900	b
0.5470	10000.0800	b
0.5365	10000.0700	b
0.6513	10000.0100	b
0.6230	10000.0000	b
0.4730	10000.0300	b
0.5084	10000.0300	b
0.6013	10000.0900	b
0.6687	10000.1200	b
0.7157	10000.0300	b
0.5481	10000.1600	b
0.5372	10000.0400	b
0.4341	10000.0700	b
0.7381	10000.1100	b
0.7409	10000.0600	b
0.6652	10000.0900	b
0.6207	10000.0300	b
0.4577	10000.0500	b
0.7014	10000.2100	b
0.6515	10000.0200	b
0.6397	10000.0500	b
0.5356	10000.0300	b
0.5798	10000.0000	b
0.4970	10000.3000	b
0.8712	10000.0700	b
0.6035	10000.0100	b
0.5963	10000.0200	b
0.6082	10000.4200	b
0.4083	10000.1000	b
0.6196	10000.0900	b
0.6309	10000.0800	b
0.5674	10000.1000	b
0.6298	10000.1900	b
0.6366	10000.0800	b
0.5177	10000.1000	b
0.5584	10000.0600	b
0.2650	55.4300	b
0.2715	0.8000	b
0.2479	55.3300	b
0.4861	239.1400	b
0.2620	0.6300	b
0.3636	238.9400	b
0.4337	18.3900	b
0.2585	0.5600	b
0.4213	18.4000	b
0.3581	20.8700	b
0.2632	0.5600	b
0.3900	20.9100	b
0.3555	36.9200	b
0.2652	0.9800	b
0.3407	36.9300	b
0.4621	65.4800	b
0.2653	0.8600	b
0.5027	65.2200	b
0.3680	19.8400	b
0.3916	0.5700	b
0.4101	19.8300	b
0.5199	16.9100	b
0.2593	0.5900	b
0.4770	16.9100	b
0.5783	37.2000	b
0.2663	0.6500	b
0.5391	37.2000	b
0.6866	226.9500	b
0.2558	0.5400	b
0.4723	227.2500	b
};
\draw [dashed] (-8.5,-10) -- (95,100);
	\end{loglogaxis}
\end{tikzpicture}
            \caption[]%
            {\qIP\ vs \BARON.\vspace{0.3\baselineskip}  \quad}    
            \label{fig:SQP35_BA}
        \end{subfigure}
        \begin{subfigure}{0.495\textwidth}   
            \centering 
            \begin{tikzpicture}[scale = \myplotscale]
	\begin{loglogaxis}[%
	height = 8.5cm,
	width = 10cm,
	xmin=0,
  xmax=12000,
  ymin=0,
  ymax=12000,
    xlabel={{\tt quadprogIP}},
    ylabel={{\tt CPLEX}},
	scatter/classes={%
		a={mark=square,black, scale = 0.75},
		b={mark=square,black, scale = 1.2},
		c={mark=square,black, scale = 1.8}}]
	\addplot[scatter,only marks,%
		scatter src=explicit symbolic]%
	table[meta=label] {
x     y      label
0.4975	46.4115	a
0.3817	46.2745	a
0.4375	10634.4272	a
0.4396	102.7590	a
0.4778	102.0703	a
0.6686	10816.7864	a
0.6939	38.4305	a
0.4941	39.3925	a
0.6090	10765.4349	a
0.4263	31.8153	a
0.4109	32.2296	a
0.6045	10889.0082	a
0.6384	16.1331	a
0.6116	15.6528	a
0.6055	10719.7743	a
0.6860	136.1634	a
0.7164	134.8724	a
0.7044	6367.2862	a
0.4895	39.5464	a
0.4960	39.6472	a
0.6905	10845.6697	a
0.5285	20.5130	a
0.4506	19.1355	a
0.4459	10848.6388	a
0.4381	46.9095	a
0.4150	46.6656	a
0.9468	10949.3666	a
0.4784	58.7189	a
0.4409	58.8939	a
0.6390	10666.1261	a
0.4146	4.5954	a
0.5185	1.4745	a
0.4691	8.6734	a
0.3999	1.9123	a
0.4415	2.2058	a
0.5341	3.3069	a
0.3426	4.2356	a
0.4301	1.4054	a
0.3838	5.9181	a
0.5270	1.4168	a
0.4854	1.3211	a
0.4447	1.7470	a
0.4456	4.1347	a
0.4147	2.1178	a
0.4710	8.7865	a
0.4316	2.9620	a
0.3194	3.5433	a
0.3473	9.4870	a
0.2818	2.2758	a
0.3804	1.8342	a
0.6049	4.3228	a
0.5007	2.1381	a
0.4094	1.9210	a
0.5149	40.6701	a
0.3279	1.7711	a
0.3022	1.0621	a
0.3861	11.6896	a
0.4358	1.3556	a
0.4189	1.1928	a
0.3047	3.7783	a
0.4909	1.6686	a
0.3426	1.5235	a
0.5152	6.0488	a
0.3784	3.5945	a
0.2976	2.0788	a
0.5979	4.9979	a
0.3595	1.4339	a
0.4454	1.2692	a
0.4638	2.6326	a
0.4499	1.7017	a
0.5686	2.2206	a
0.6415	6.2288	a
0.3709	1.0336	a
0.3317	0.7333	a
0.3193	1.8207	a
0.3890	2.5692	a
0.3297	0.8637	a
0.5410	3.0070	a
0.2989	1.8935	a
0.3864	1.3829	a
0.5938	4.3347	a
0.3289	2.2095	a
0.3130	1.5491	a
0.3943	4.0475	a
0.3256	1.0415	a
0.3734	1.0543	a
0.4095	4.5168	a
0.5136	1.9932	a
0.3483	1.1802	a
0.4060	1.9723	a
0.3052	1.1110	a
0.2994	1.2124	a
0.3193	1.9121	a
0.3336	1.5849	a
0.2869	0.8908	a
0.3745	1.6749	a
0.3105	2.9841	a
0.6704	0.8027	a
0.3138	3.5343	a
0.3090	1.3850	a
0.5048	0.7948	a
0.3341	1.3976	a
0.4274	1.1516	a
0.3268	0.9166	a
0.4323	3.0218	a
0.5004	1.2376	a
0.4342	0.7179	a
0.4084	1.3962	a
0.4338	1.0179	a
0.4088	0.8148	a
0.3662	0.9017	a
0.3180	1.9475	a
0.3826	1.0986	a
0.5515	1.8831	a
0.3351	1.1803	a
0.3175	1.2681	a
0.3629	1.4899	a
0.4690	1.1238	a
0.4677	0.9341	a
0.2791	1.7454	a
0.2388	1.7740	a
0.2414	0.2160	a
0.2325	1.4913	a
0.4004	9.4156	a
0.2256	0.2912	a
0.4837	9.3485	a
0.2365	1.5039	a
0.2310	0.2498	a
0.2483	1.4074	a
0.3962	4.4515	a
0.2413	0.2625	a
0.3240	4.6847	a
0.3677	4.8787	a
0.2294	0.1738	a
0.4230	5.0948	a
0.4825	19.9287	a
0.2420	0.3181	a
0.2919	19.9189	a
0.4771	1.6551	a
0.2319	0.1325	a
0.5767	1.5749	a
0.3005	4.3887	a
0.2357	0.1518	a
0.3288	4.4584	a
0.2309	0.4906	a
0.2320	0.1731	a
0.2247	0.4198	a
0.3118	3.8887	a
0.2253	0.2490	a
0.3588	3.8026	a
		
0.9329	5559.7759	b
0.8337	5561.6415	b
4.1844	10130.3488	b
0.6683	2134.4316	b
0.6790	2154.5814	b
4.1127	10114.5317	b
0.9613	426.5903	b
0.9209	421.9571	b
3.0645	10121.5475	b
0.8144	2505.8531	b
0.7441	2505.8643	b
3.3178	10145.0042	b
0.8256	8860.0751	b
0.7794	8911.7344	b
4.4178	10084.0362	b
0.6432	4104.0470	b
0.5462	4116.7452	b
3.7384	10119.0365	b
0.7037	200.0723	b
0.7170	199.1931	b
3.1523	10114.9630	b
1.2164	4219.8061	b
1.2485	8115.1763	b
31.1552	10072.9925	b
0.8713	9126.4004	b
13.4146	10034.4355	b
14.9197	10079.4140	b
0.5918	4825.9605	b
0.7356	8760.6838	b
25.9848	10068.0135	b
0.5340	14.7408	b
0.4001	27.4874	b
0.8275	290.4149	b
0.5547	14.3508	b
0.4433	33.0412	b
0.5355	183.5355	b
0.6157	9.1945	b
0.5353	9.3775	b
0.5723	326.0799	b
0.4247	11.6858	b
0.8469	8.5111	b
0.4950	45.7871	b
0.5373	59.9680	b
0.4464	18.3502	b
0.5460	170.6784	b
0.5244	33.8226	b
0.4621	39.6354	b
0.8702	300.2541	b
0.9726	47.3054	b
0.6979	61.5629	b
0.5182	667.4581	b
0.5292	28.9258	b
0.4901	23.5810	b
0.6915	109.1766	b
0.4844	76.7392	b
0.6172	71.8251	b
0.8885	234.2369	b
0.8747	45.8403	b
0.7515	14.7879	b
0.6005	397.6179	b
0.5450	7.7981	b
0.6686	17.8243	b
0.5442	36.8422	b
0.6333	15.0771	b
0.5831	5.0592	b
0.6157	53.3695	b
0.4218	5.4080	b
0.5242	4.8892	b
0.4514	61.9289	b
0.5321	15.0726	b
0.5512	16.8395	b
0.5921	67.9708	b
0.5857	31.3586	b
0.5615	10.8836	b
0.5537	71.2410	b
0.4100	40.3045	b
0.6200	6.6290	b
0.4574	40.9897	b
0.5824	15.5641	b
0.5834	8.4140	b
0.4646	51.8731	b
0.5583	9.9713	b
0.7534	5.2296	b
0.5359	52.6987	b
0.8077	22.6528	b
0.5470	7.9890	b
0.5365	97.3018	b
0.6513	7.9964	b
0.6230	6.7831	b
0.4730	45.4681	b
0.5084	14.3318	b
0.6013	5.7499	b
0.6687	29.5028	b
0.7157	6.3187	b
0.5481	7.7556	b
0.5372	35.6256	b
0.4341	8.1102	b
0.7381	4.7415	b
0.7409	31.1481	b
0.6652	11.5061	b
0.6207	6.9633	b
0.4577	26.5755	b
0.7014	4.5871	b
0.6515	5.4341	b
0.6397	36.3641	b
0.5356	8.3110	b
0.5798	4.3796	b
0.4970	9.4189	b
0.8712	5.6128	b
0.6035	3.0458	b
0.5963	35.5189	b
0.6082	4.3992	b
0.4083	3.2226	b
0.6196	8.6940	b
0.6309	15.2571	b
0.5674	4.2192	b
0.6298	7.2481	b
0.6366	14.6039	b
0.5177	6.1000	b
0.5584	48.7162	b
0.2650	7.2497	b
0.2715	0.2582	b
0.2479	7.2247	b
0.4861	17.8437	b
0.2620	0.2304	b
0.3636	18.0479	b
0.4337	18.2657	b
0.2585	0.0262	b
0.4213	18.5905	b
0.3581	9.1112	b
0.2632	0.2815	b
0.3900	9.0962	b
0.3555	75.2451	b
0.2652	0.1134	b
0.3407	79.7035	b
0.4621	84.9314	b
0.2653	0.3656	b
0.5027	84.0155	b
0.3680	2.7064	b
0.3916	0.2837	b
0.4101	2.6402	b
0.5199	3.1382	b
0.2593	0.3835	b
0.4770	3.0650	b
0.5783	100.7327	b
0.2663	0.1312	b
0.5391	96.0623	b
0.6866	25.9695	b
0.2558	0.2532	b
0.4723	25.4609	b
};
\draw [dashed] (-10,-10) -- (120,120);
	\end{loglogaxis}
\end{tikzpicture}
            \caption[]%
            {\qIP\ vs \cplex.\vspace{0.3\baselineskip}  \quad}    
            \label{fig:SQP35_CP}
        \end{subfigure}
        \begin{subfigure}{0.495\textwidth}   
            \centering 
            \includegraphics[width=\textwidth]{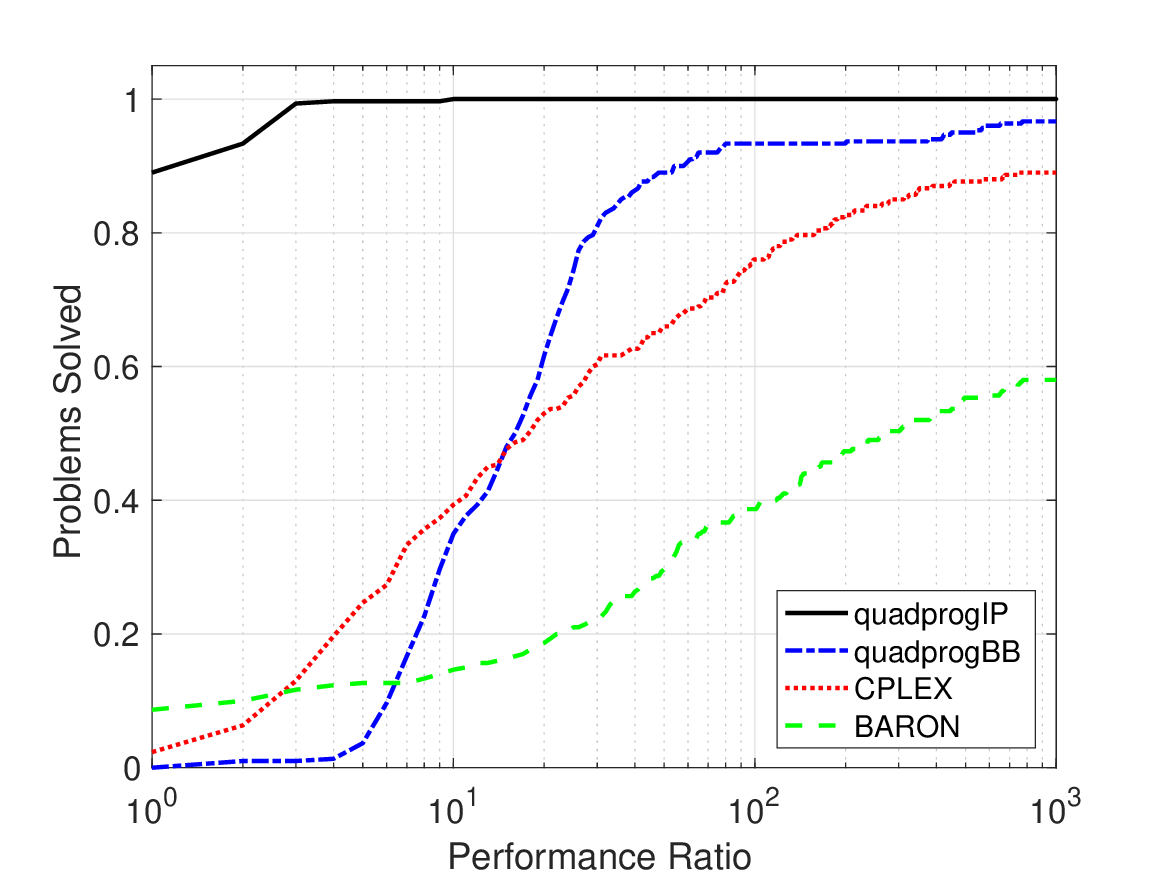}
            \caption[]%
            {Performance profile for SQP30 and SQP50 instances.\vspace{0.3\baselineskip}  \quad}    
            \label{fig:SQP35_Profile}
        \end{subfigure}
        \caption[]
        {Solution time in seconds of SQP30 and SQP50 instances. Size of squares illustrates size of the instance. A square at the maximum value of an axis represents an instances for which the solver in that axis reached maximum running time without solving it.} 
        \label{fig:SQP35_ALL}
    \end{figure}

As Figure~\ref{fig:SQP35_ALL} shows, the results on the SQP instances SQP30 and SQP50 from~\citet{bonami2016solving} is very similar to the ones presented in Section~\ref{sec:SQPnum}. As with the set of SQP instances, only \cplex is able to solve a few instances faster than \qIP; however, in general \qIP outperforms the other solvers by orders of magnitude in terms of solution time.

\subsubsection{Results on StableQP instances.}

In line with the performance of \qIP\  on SQP, SQP30, and SQP50 instances, it is
interesting to see in Table~\ref{tab:stable} that \qIP\  clearly outperforms \qBB, \BARON, and \cplex in the StableQP instances (see, Section~\ref{sec:instances}).
In fact, while \qIP\  solves each of the instances in less than a second, \qBB, and \cplex are unable to solve the instances beyond $k \ge 4$ within
the maximum allowed solution time of $10^4$ seconds, while \BARON{} is unable to solve the instances beyond $k\ge 3$ within the maximum allowed solution time.

\begin{table}[htb]
\begin{center}
{\small
  \begin{tabular}{ c  r r r r}
  \toprule
  & \multicolumn{4}{c}{Solution Time (s.)}\\
  \cline{2-5}
     $k$ & \multicolumn{1}{c}{\tt quadprogIP} & \multicolumn{1}{c}{\tt quadprogBB} & \multicolumn{1}{c}{\tt\t  \BARON{}} & \multicolumn{1}{c}{\tt\t  \cplex} \\
     \midrule
         1 & 0.34 & 3.67 & 8.93 & 0.39\\
    2 & 0.25 & 6.28 & 2573.77 & 8.75\\
        3 & 0.34 & 12.56 & - & 685.70\\
    4 &  0.43& - & - & -\\
    5 &  0.49&   - & - & -\\
    6 &  0.51&  -  & - & -\\
    7 &  0.46&   - &-  & -\\
    8 &  0.49&   - & -  & -\\
    \bottomrule
  \end{tabular}
  }
  \end{center}
  \caption{Solution time in seconds for StableQP instances. Dash ``-'' indicates that the solver was unable to solve the instance within the maximum allowed time.}
  \label{tab:stable}
\end{table}

\subsubsection{Results on Scozzari/Tardella instances.} The Scozzari/Tardella from~\citet{scozzari2008clique} are composed of much larger-scale instances of SQP than the ones considered so far. Table~\ref{tab:tardella}, as with the previously discussed groups of standard QP instances, clearly shows that \qIP is able to solve these instances faster than the other solvers, and is able to solve more large-scale instances than the other solvers.

\begin{table}[!htb]
\begin{center}
{\small
  \begin{tabular}{ l  r r r r}
  \toprule
   & \multicolumn{4}{c}{Solution Time (s.)}\\
  \cline{2-5}
     Scozzari/Tardella instance & \multicolumn{1}{c}{\tt quadprogIP} & \multicolumn{1}{c}{\tt quadprogBB} & \multicolumn{1}{c}{\tt\t  \BARON{} } & \multicolumn{1}{c}{\tt\t \cplex} \\
     \midrule
Problem\_30x30\_0.75.mps.mat:	&	 0.51 	&	 5.27 	&	 39.32 	&	  5.56	\\     
Problem\_50x50\_0.75.mps.mat:	&	 11.78 	&	 48.29 	&	 - 	&	 2,162.13 	\\
Problem\_100x100\_0-1.mps.mat:	&	 1.54 	&	 1,412.16 	&	 223.54 	&	 154.71 	\\
Problem\_100x100\_0.5.mps.mat:	&	 6.71 	&	 319.82 	&	 - 	&	 - 	\\
Problem\_100x100\_0.75.mps.mat:	&	 36.76 	&	 1,519.61 	&	 - 	&	 - 	\\
Problem\_200x200\_0-1.mps.mat:	&	 36.86 	&	 - 	&	 - 	&	 9,995.67 	\\
Problem\_200x200\_0.5.mps.mat:	&	 175.15 	&	 - 	&	 - 	&	 - 	\\
Problem\_500x500\_0-1.mps.mat:	&	 240.09 	&	 - 	&	 - 	&	 - 	\\
Problem\_500x500\_0.25.mps.mat:	&	 2,092.48 	&	 - 	&	 - 	&	 - 	\\
Problem\_1000x1000\_0.25.mps.mat:	&	 - 	&	 - 	&	 - 	&	 - 	\\
Problem\_Q30.mps.mat:	&	 0.54 	&	 4.27 	&	 - 	&	 - 	\\
Problem\_Q50.mps.bar.mat:	&	 2.45 	&	 8,476.70 	&	 - 	&	 - 	\\
Problem\_Q100.mps.bar.mat:	&	 4.71 	&	 - 	&	 - 	&	 - 	\\
Problem\_Q150.mps.mat:	&	 20.89 	&	 - 	&	 - 	&	 - 	\\    \bottomrule
  \end{tabular}
  }
  \end{center}
  \caption{Solution time in seconds for Scozzari/Tardella instances. Dash ``-'' indicates that the solver was unable to solve the instance within the maximum allowed time.}
  \label{tab:tardella}
\end{table}

\subsubsection{Results on BoxQP instances.} 
\begin{figure}[!htb]
    \captionsetup[figure]{justification=left}
    \captionsetup[subfigure]{justification=centering}
        \centering
        \begin{subfigure}{0.495\textwidth}
            \centering
            \begin{tikzpicture}[scale = \myplotscale]
	\begin{loglogaxis}[%
	height = 8.5cm,
	width = 10cm,
	xmin=0,
  xmax=12000,
  ymin=0,
  ymax=12000,
    xlabel={{\tt quadprogIP}},
    ylabel={{\tt quadprogBB}},
	scatter/classes={%
		a={mark=square,black, scale = 0.75},
		b={mark=square,black, scale = 1.2},
		c={mark=square,black, scale = 1.8}}]
	\addplot[scatter,only marks,%
		scatter src=explicit symbolic]%
	table[meta=label] {
x     y      label
0.544653	8.047237	a
0.471305	4.941475	a
0.520606	1.951007	a
10.644904	24.480365	a
0.723876	5.599081	a
2.737187	28.06778	a
27.627206	24.110175	a
1.256076	8.628496	a
1.384003	39.307513	a
20.288066	29.976201	a
0.696845	3.74154	a
0.564317	5.747346	a
2.156107	4.024761	a
2.405927	18.2434	a
1.415431	3.894412	a
12.147442	5.872524	a
18.863345	8.908819	a
2.81134	25.925627	a
1.502259	17.130322	b
1.469603	18.682278	b
1.698555	23.947138	b
70.745784	87.949358	b
2.878179	7.65249	b
50.442611	30.511641	b
50.86874	48.053503	b
21.615357	76.756976	b
24.255243	35.098249	b
279.897035	265.207718	b
4.622095	28.498149	b
2.553951	7.081098	b
39.641255	17.349613	b
38.107733	9.159534	b
22.404305	23.71587	b
136.497692	8.355219	b
97.980438	9.577865	b
23.600968	26.008089	b
90.548894	20.965029	b
114.861873	17.12126	b
41.395776	8.529093	b
50.824202	20.4322	b
433.762402	36.398354	b
10002.55075	111.025775	b
6.350487	13.964336	b
74.568404	60.814007	b
51.58166	101.198434	b
219.587484	28.293257	b
363.38801	99.581398	b
67.183777	25.726996	b
10020.14194	844.37019	b
9362.434466	186.915374	b
10002.81242	151.018886	b
73.568063	61.97576	b
25.946898	120.076665	b
307.0833	122.559369	b
10007.18145	303.652691	c
10012.37114	225.248516	c
10007.48821	433.582149	c
10013.77368	303.669808	c
10012.53704	176.068405	c
10008.05911	51.102252	c
10014.10815	421.355301	c
10014.20897	1975.019126	c
10013.3445	930.342406	c
10010.57053	194.124165	c
10013.47276	561.830679	c
10010.05485	649.853316	c
10012.00733	7036.097202	c
10011.59906	222.220238	c
10011.44875	961.183402	c
10012.16985	1020.159782	c
10012.44875	1593.863007	c
10012.84633	4039.243716	c
10013.04942	1873.956998	c
10008.97597	1437.047613	c
10009.97221	1144.723442	c
10012.42638	3306.96082	c
10009.76267	702.288253	c
10011.258	2308.840278	c
10011.79919	9078.976486	c
10011.16903	8842.309739	c
10011.17656	3562.650047	c
10008.57676	1507.059432	c
10008.27423	1103.152197	c
10009.62271	1253.511631	c
10009.54751	10003.80579	c
10009.71666	10000.47822	c
10010.43318	2153.699966	c
10010.6643	10007.02291	c
10009.41216	10001.82576	c
10012.68021	10003.57664	c
};
\draw [dashed] (0,0) -- (120,120);
	\end{loglogaxis}
\end{tikzpicture}
            \caption{\qIP\ vs \qBB.\vspace{0.3\baselineskip} \quad}    
            \label{fig:BoxQP_BB}
        \end{subfigure}
        \begin{subfigure}{0.495\textwidth}  
            \centering 
            \begin{tikzpicture}[scale = \myplotscale]
	\begin{loglogaxis}[%
	height = 8.5cm,
	width = 10cm,
	xmin=0,
  xmax=12000,
  ymin=0,
  ymax=12000,
    xlabel={{\tt quadprogIP}},
    ylabel={{\tt \BARON{}}},
	scatter/classes={%
		a={mark=square,black, scale = 0.75},
		b={mark=square,black, scale = 1.2},
		c={mark=square,black, scale = 1.8}}]
	\addplot[scatter,only marks,%
		scatter src=explicit symbolic]%
	table[meta=label] {
x     y      label
0.544653	0.26	a
0.471305	0.9	a
0.520606	0.41	a
10.644904	1.32	a
0.723876	0.66	a
2.737187	1.91	a
27.627206	7.2	a
1.256076	0.75	a
1.384003	2.16	a
20.288066	4.43	a
0.696845	0.39	a
0.564317	1.21	a
2.156107	3.92	a
2.405927	3.74	a
1.415431	0.94	a
12.147442	11.7	a
18.863345	3.03	a
2.81134	16.18	a
1.502259	0.24	b
1.469603	0.85	b
1.698555	0.42	b
70.745784	3.65	b
2.878179	0.73	b
50.442611	1.97	b
50.86874	6.8	b
21.615357	2.05	b
24.255243	4.62	b
279.897035	47.23	b
4.622095	3.3	b
2.553951	0.61	b
39.641255	22.32	b
38.107733	16.6	b
22.404305	7.67	b
136.497692	24.85	b
97.980438	14.29	b
23.600968	10.02	b
90.548894	33.17	b
114.861873	3.38	b
41.395776	20.13	b
50.824202	5.25	b
433.762402	135.68	b
10002.55075	3464.9	b
6.350487	0.46	b
74.568404	0.98	b
51.58166	1.24	b
219.587484	3.48	b
363.38801	5	b
67.183777	6.12	b
10020.14194	151.53	b
9362.434466	9.2	b
10002.81242	14.81	b
73.568063	1.23	b
25.946898	0.6	b
307.0833	0.74	b
10007.18145	12.1	c
10012.37114	25.2	c
10007.48821	11.85	c
10013.77368	596.89	c
10012.53704	153.62	c
10008.05911	6.8	c
10014.10815	10000.02	c
10014.20897	10000.17	c
10013.3445	10000.19	c
10010.57053	13.31	c
10013.47276	651.05	c
10010.05485	47.6	c
10012.00733	10000.01	c
10011.59906	419.26	c
10011.44875	333.29	c
10012.16985	10000.33	c
10012.44875	10000.31	c
10012.84633	10000.25	c
10013.04942	192	c
10008.97597	212.72	c
10009.97221	131.82	c
10012.42638	10000.17	c
10009.76267	795.44	c
10011.258	9413.98	c
10011.79919	10000.09	c
10011.16903	10000.47	c
10011.17656	10000.33	c
10008.57676	3713.21	c
10008.27423	988.46	c
10009.62271	676.99	c
10009.54751	10000.2	c
10009.71666	10000.23	c
10010.43318	10000.33	c
10010.6643	10000.49	c
10009.41216	10000.46	c
10012.68021	10000.01	c
};
\draw [dashed] (-8.5,-10) -- (95,100);
	\end{loglogaxis}
\end{tikzpicture}
            \caption[]%
            {\qIP\ vs \BARON.\vspace{0.3\baselineskip}  \quad}    
            \label{fig:BoxQP_BA}
        \end{subfigure}
        \begin{subfigure}{0.495\textwidth}   
            \centering 
            \begin{tikzpicture}[scale = \myplotscale]
	\begin{loglogaxis}[%
	height = 8.5cm,
	width = 10cm,
	xmin=0,
  xmax=12000,
  ymin=0,
  ymax=12000,
    xlabel={{\tt quadprogIP}},
    ylabel={{\tt CPLEX}},
	scatter/classes={%
		a={mark=square,black, scale = 0.75},
		b={mark=square,black, scale = 1.2},
		c={mark=square,black, scale = 1.8}}]
	\addplot[scatter,only marks,%
		scatter src=explicit symbolic]%
	table[meta=label] {
x     y      label
0.544653	0.137683	a
0.471305	0.358666	a
0.520606	0.391395	a
10.644904	1.3774	a
0.723876	0.365303	a
2.737187	0.943832	a
27.627206	1.96684	a
1.256076	0.534089	a
1.384003	0.562099	a
20.288066	2.831874	a
0.696845	0.179849	a
0.564317	0.208522	a
2.156107	1.018643	a
2.405927	1.233669	a
1.415431	0.609508	a
12.147442	1.806741	a
18.863345	2.930743	a
2.81134	1.739948	a
1.502259	0.058206	b
1.469603	0.253923	b
1.698555	0.150644	b
70.745784	2.287916	b
2.878179	0.187658	b
50.442611	2.149615	b
50.86874	1.624174	b
21.615357	0.898955	b
24.255243	1.011193	b
279.897035	4.32951	b
4.622095	0.955889	b
2.553951	0.756813	b
39.641255	3.182198	b
38.107733	1.696757	b
22.404305	1.37421	b
136.497692	4.082391	b
97.980438	2.917697	b
23.600968	3.858319	b
90.548894	4.776801	b
114.861873	8.504264	b
41.395776	3.543792	b
50.824202	5.757698	b
433.762402	13.227284	b
10002.55075	143.336346	b
6.350487	0.225482	b
74.568404	1.668153	b
51.58166	1.165925	b
219.587484	3.890638	b
363.38801	3.01222	b
67.183777	1.512964	b
10020.14194	59.036939	b
9362.434466	6.262719	b
10002.81242	8.565442	b
73.568063	0.452399	b
25.946898	0.321514	b
307.0833	1.35418	b
10007.18145	2.986466	c
10012.37114	6.279143	c
10007.48821	7.003652	c
10013.77368	75.616337	c
10012.53704	36.064413	c
10008.05911	11.149074	c
10014.10815	583.271432	c
10014.20897	10000.08646	c
10013.3445	5527.538183	c
10010.57053	5.043635	c
10013.47276	21.179438	c
10010.05485	12.371608	c
10012.00733	10000.25583	c
10011.59906	79.978732	c
10011.44875	125.123883	c
10012.16985	10000.05374	c
10012.44875	10000.05848	c
10012.84633	10000.06774	c
10013.04942	47.463502	c
10008.97597	40.325211	c
10009.97221	27.981375	c
10012.42638	2653.659052	c
10009.76267	250.203434	c
10011.258	261.0491	c
10011.79919	10000.02866	c
10011.16903	10000.05362	c
10011.17656	10000.03891	c
10008.57676	108.76555	c
10008.27423	74.77358	c
10009.62271	64.969638	c
10009.54751	10000.04447	c
10009.71666	10000.07744	c
10010.43318	10000.06439	c
10010.6643	10000.06595	c
10009.41216	10000.65648	c
10012.68021	10000.02963	c
};
\draw [dashed] (-10,-10) -- (120,120);
	\end{loglogaxis}
\end{tikzpicture}
            \caption[]%
            {\qIP\ vs \cplex.\vspace{0.3\baselineskip}  \quad}    
            \label{fig:BoxQP_CP}
        \end{subfigure}
        \begin{subfigure}{0.495\textwidth}   
            \centering 
            \includegraphics[width=\textwidth]{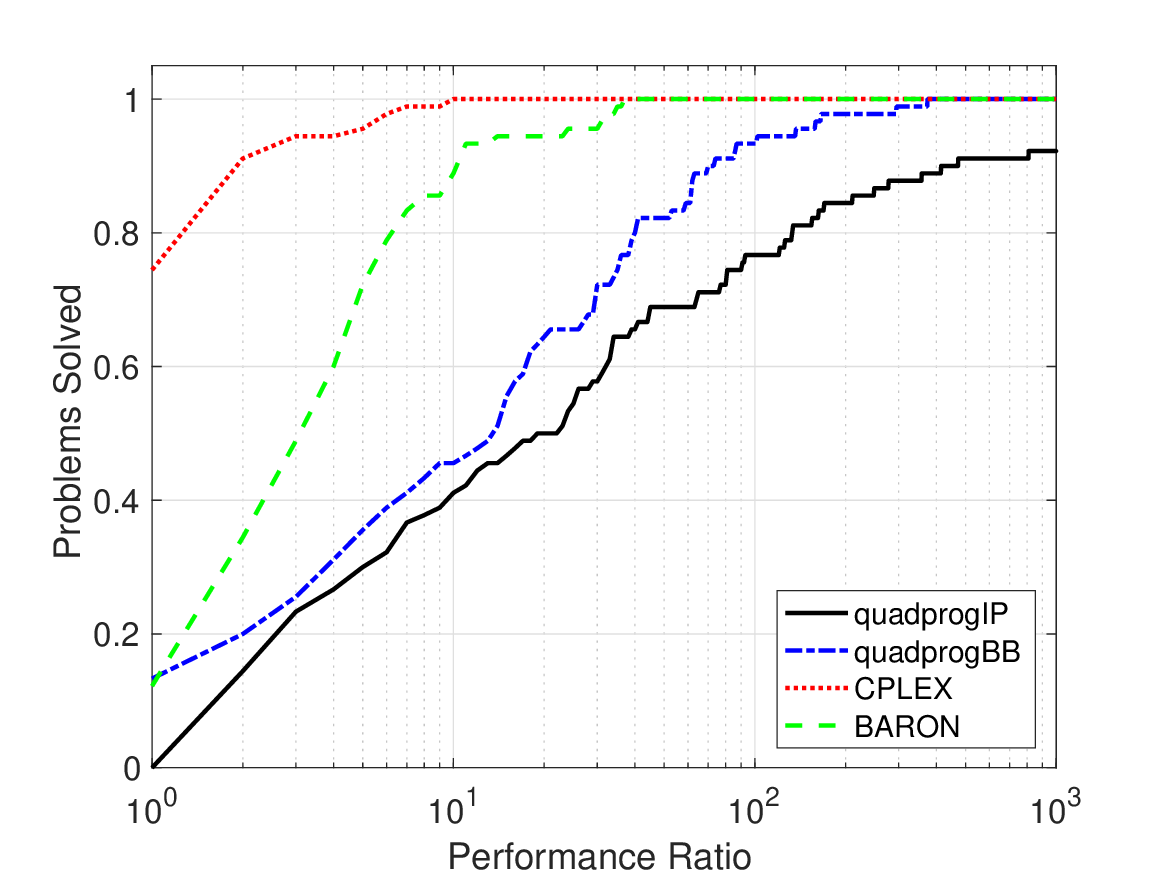}
            \caption[]%
            {Performance profile for BoxQP instances.\vspace{0.3\baselineskip}  \quad}    
            \label{fig:BoxQP_Profile}
        \end{subfigure}
        \caption[]
        {Solution time in seconds of BoxQP instances. Size of squares illustrates size of the instance. A square at the maximum value of an axis represents an instances for which the solver in that axis reached maximum running time without solving it.} 
        \label{fig:BoxQP_ALL}
    \end{figure}
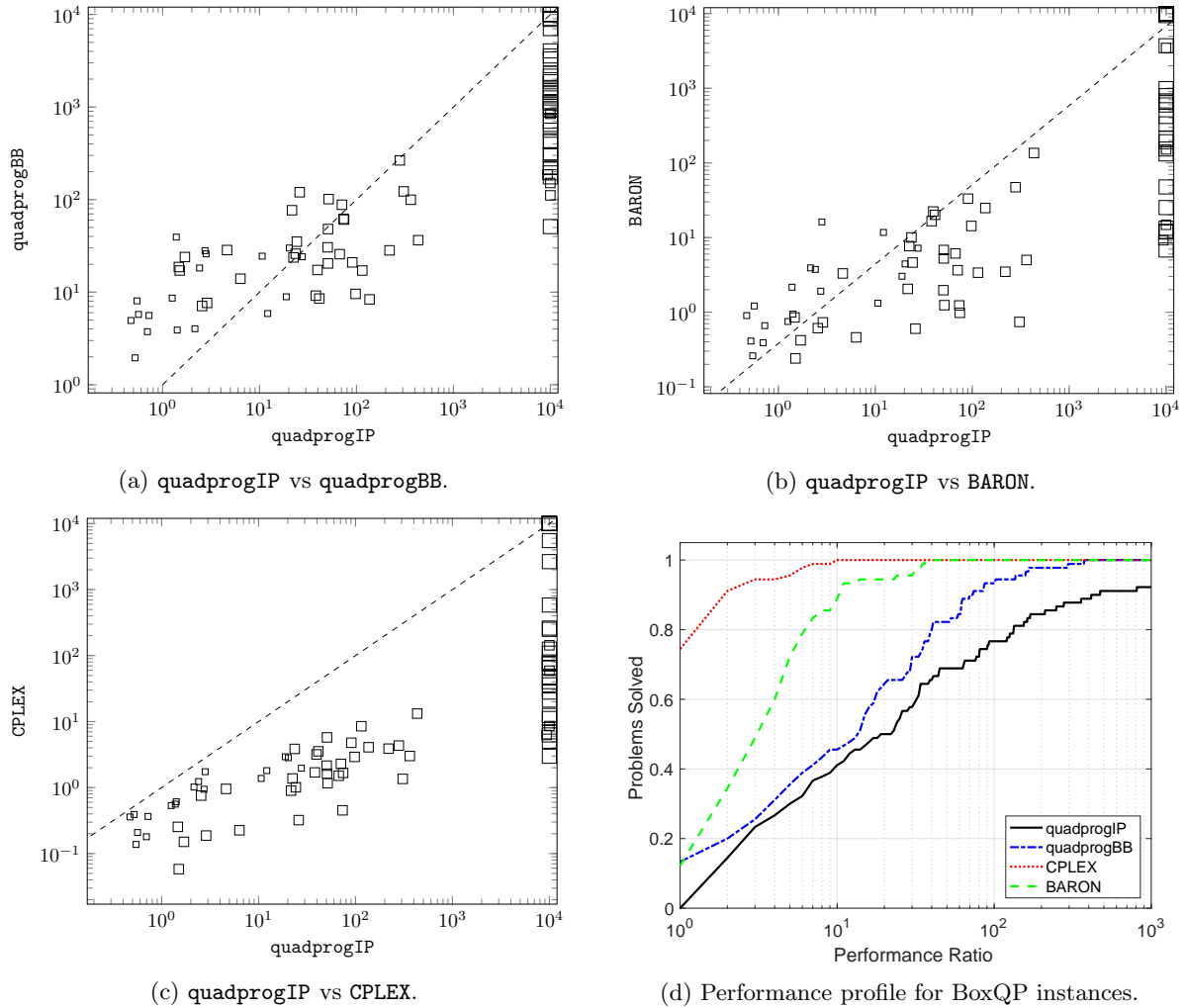

In Figure~\ref{fig:BoxQP_ALL}, we compare the performance of \qIP\  on the BoxQP instances against the other three selected solvers. It is clear from
Figure~\ref{fig:BoxQP_BB} that while \qIP\  outperforms \qBB\ in the smaller BoxQP instances (ranging between 20--60
decision variables), \qBB\ outperforms \qIP\  for larger BoxQP instances (ranging between
60--100 decision variables), where \qIP\  is typically unable to solve the instance within
the $10^4$ maximum solution time.

Figure~\ref{fig:BoxQP_BA} shows the performance of \qIP\ and \BARON{} on the BoxQP test set. It is clear that \BARON{} outperforms \qIP\  in most BoxQP instances. Although for instances with less than 40 decision variables the solution time of \qIP\ is not significantly longer than that of \BARON{}. Figure~\ref{fig:BoxQP_CP} shows that \cplex performs much better than \qIP on all BoxQP instances. Figure~\ref{fig:BoxQP_Profile} summarizes these results, where it is clear that \cplex and \BARON are the best solvers for these BoxQP instances.

It is worth mentioning that the performance of \qIP on BoxQP instances can be improved by adding appropriate valid constraints to the IQP~\eqref{eq:milpQP} formulation of the BoxQP.
This valid constraints can be derived from~\citet[][Prop. 1]{hansen1993global}. Specifically, notice that the~IQP~\eqref{eq:milpQP} corresponding to~\eqref{eq:BQP_std}
can be written as:

\begin{equation}
\label{BoxmilpQP}
\begin{array}{rllll}
\frac 12\min &  (Hl+f)\tr  x  - (u-l)\tr  \mu \\
\st &  \begin{bmatrix}H & 0\\0 & 0\end{bmatrix}\begin{bmatrix}x\\s\end{bmatrix}+\begin{bmatrix}f +Hl\\0\end{bmatrix}+\begin{bmatrix}I \\ I\end{bmatrix}\mu - \begin{bmatrix}\lambda^x\\\lambda^s\end{bmatrix} = 0\\
&x + s = u-l \\
&  0 \le x_i \le z^x_i (u_i-l_i) & i=1,\dots,n\\
&  0 \le s_i \le z^s_i (u_i-l_i) & i=1,\dots,n\\
&  0 \le \lambda^x_i \le (1-z^x_i) V_i & i=1,\dots,n\\
&  0 \le \lambda^s_i \le (1-z^s_i) V_i & i=1,\dots,n\\
&  z^x_i, z^s_i \in \{0,1\}& i=1,\dots,n.
\end{array}
\end{equation}
Then, from~\citet[][Prop. 1]{hansen1993global}, and Proposition~\ref{prop.boxbound} below, it follows that the constraints
\begin{equation}
\label{eq.valid}
z^x_i + z^s_i = 1, \text{for all $i \in \{1,\dots,n\}$ such that $H_{ii} < 0$,}
\end{equation}
are valid constraints for the optimal solutions of~\eqref{BoxmilpQP}. When added to~\eqref{BoxmilpQP}, the valid constraints~\eqref{eq.valid}
improve the solution time of the approach proposed here to globally solve BoxQP problems.

\begin{proposition}
\label{prop.boxbound}
Let $(x, s, \mu, \lambda^x, \lambda^s) \in \R^n_+ \times \R^n_+ \times \R^{2n} \times \R^n_+ \times \R^n_+$ be
a KKT point of~\eqref{eq:BQP_std}. Then  $e\tr \lambda^x + e\tr \lambda^s \le M$, where $M$ is given by~\eqref{eq:boxbd}. That is, every KKT point of \eqref{eq:BQP_std} satisfies the bounds given by~\eqref{eq:boxbd}.
\end{proposition}

\proof{Proof.}
Let $(x, s, \mu, \lambda^x, \lambda^s)$ satisfy the KKT conditions of~\eqref{eq:BQP_std}. Then
\begin{align}
\begin{bmatrix}H & 0\\0 & 0\end{bmatrix}\begin{bmatrix}x\\s\end{bmatrix}+\begin{bmatrix}f +Hl\\0\end{bmatrix}+\begin{bmatrix}I \\ I\end{bmatrix}\mu - \begin{bmatrix}\lambda^x\\\lambda^s\end{bmatrix} = 0\label{KKT}\\
x\tr \lambda^x = 0, s\tr \lambda^s = 0,
x+s = U \label{comp_feas}\\ x\geq 0, \lambda^x \geq 0, \lambda^s \geq 0 \notag.
\end{align}
where $\mu \in \R^n, \lambda^x\in \mathbb{R}^n_+, \lambda^x\in \mathbb{R}^n_+$ are respectively the dual multipliers for the $x+s = u-l$, $x \ge 0$,
and $s \ge 0$ constraints in~\eqref{eq:BQP_std}. Equation \eqref{KKT} implies that
\begin{equation}
\label{eq.lambdabox}
\lambda^x - \lambda^s = Hx + f + Hl.
\end{equation}
 Thus, using \eqref{comp_feas}, \eqref{eq.lambdabox}, and the results in Section~\ref{sec:boxqp}
\begin{align*}
e\tr \lambda^x + e\tr \lambda^s  & = e\tr |\lambda^x - \lambda^s|  =  \|Hx + f + Hl \|_1 \le  \|Hx\|_1 + \|f + Hl \|_1 \le M
 \end{align*}
\endproof{}

\begin{figure}[!htb]  
\centering 
 \includegraphics[width=0.5\textwidth]{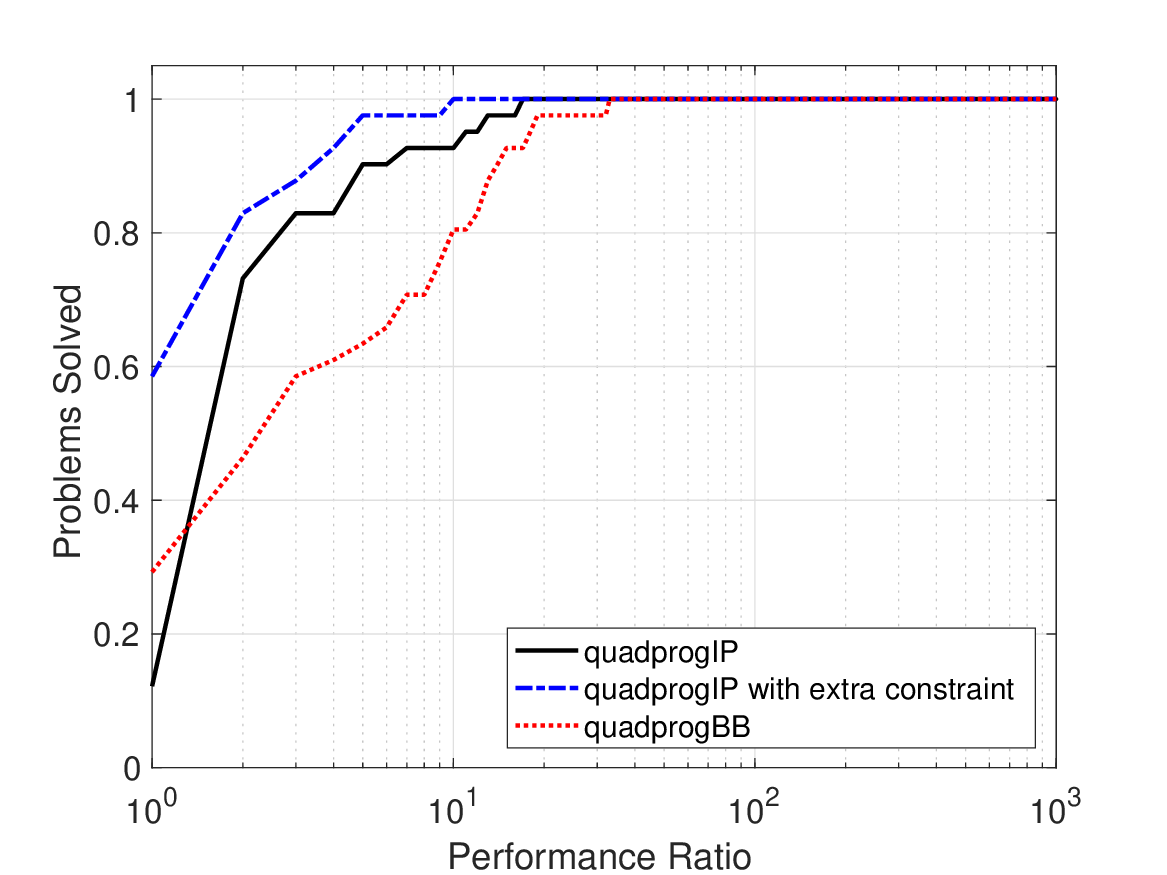}
\caption[]%
 {Performance profile for {\tt spar} BoxQP instances anging on size between 20-40 variables with density between 30-100. Extra constraints refer to adding constraints~\eqref{eq.valid} in the \qIP solver.\vspace{0.3\baselineskip}  \quad}    
 \label{fig:withcons}
   \end{figure}

Although the \qIP code does not include the strengthening constraints~\eqref{eq.valid} for BoxQPs, the results illustrated on Figure~\ref{fig:withcons} show how adding the valid constraints~\eqref{eq.valid} improves the solution time on a set of {\tt spar} BoxQP instances ranging on size between 20-40 variables with density between 30-100. In particular, with the addition of these constraints, \qIP outperforms \qBB on these instances.

\subsubsection*{Results on CUTEr, Globallib, and RandQP instances.}
\begin{figure}[!htb]
    \captionsetup[figure]{justification=left}
    \captionsetup[subfigure]{justification=centering}
        \centering
        \begin{subfigure}{0.495\textwidth}
            \centering
            \begin{tikzpicture}[scale = \myplotscale]
	\begin{loglogaxis}[%
	height = 8.5cm,
	width = 10cm,
	xmin=0,
  xmax=12000,
  ymin=0,
  ymax=12000,
    xlabel={{\tt quadprogIP}},
    ylabel={{\tt quadprogBB}},
	scatter/classes={%
		a={mark=square,black, scale = 0.75},
		b={mark=square,black, scale = 1.2},
		c={mark=square,black, scale = 1.8}}]
	\addplot[scatter,only marks,%
		scatter src=explicit symbolic]%
	table[meta=label] {
x     y      label
6.272094	5.018216	a
2.680431	20.330143	a
6.223431	2.4683	a
3.711449	1.741773	a
4.633106	1.138993	a
0.198511	1.782894	a
0.316398	13.489935	a
0.496636	5.21748	a
0.191113	0.454161	a
0.301215	2.561351	a
0.191289	0.479634	a
0.199246	2.417472	a
0.287471	21.079164	a
8.14685	180.657593	a
0.431253	641.431097	a
0.259599	1.060915	a
0.19084	0.365295	a
12.180735	29.171694	b
13.183473	28.83967	b
0.198109	9.555662	a
0.198224	8.516593	a
0.258417	1.071035	a
0.286848	1.042721	a
0.197471	0.902839	a
0.464387	0.927866	a
0.197509	0.906908	a
0.264805	6.711031	a
0.297753	0.889674	a
0.204523	0.93223	a
0.195014	0.996434	a
0.336474	1.016653	a
0.199258	0.77406	a
0.400184	12.081125	a
0.196603	0.426266	a
0.411414	2.392792	a
0.278064	2.805648	a
0.203168	2.436319	a
0.334686	21.155424	a
4.104725	38.147458	a
5.435408	33.517859	a
5.108108	35.466804	a
3.915188	35.984499	a
8.022685	181.193011	a
153.318189	10000.96979	c
0.215762	2.772654	a
0.202918	2.986826	a
0.250505	1.316989	a
0.19602	1.206877	a
0.256755	2.478111	a
0.203145	1.779073	a
0.199787	1.847671	a
0.200047	1.411032	a
0.273755	2.824961	a
0.340405	1.898209	a
0.270989	1.881249	a
0.38285	140.905689	a
0.197201	0.921329	a
0.283756	0.81141	a
0.879105	4.64604	a
12.502465	9.34192	a
0.2821	1.892234	a
0.323705	12.235241	a
0.295503	1.178332	a
0.192499	0.3242	a
0.588008	4.437416	a
0.382033	4.575776	a
0.295654	3.82684	a
0.243686	1.304811	a
0.267263	0.9662	a
0.268503	1.204477	a
0.31459	0.954074	a
0.27726	0.993171	a
0.256017	0.807069	a
0.250592	0.685762	a
0.309	0.887812	a
0.292584	1.831079	a
0.252223	1.84986	a
0.191352	0.647402	a
0.191123	0.64701	a
0.30116	0.755204	a
0.337077	28.847004	a
0.424575	103.808738	a
0.283947	1.502424	a
0.558561	2.647621	a
4.675576	44.24311	a
17.060661	57.508205	a
36.991265	1703.078208	b
125.205971	5091.610782	b
2.99288	0.898614	a
4.536537	0.856226	a
0.610013	415.545891	a
0.646539	103.565789	a
1.540013	152.003967	a
8.264153	1158.683397	a
0.389049	227.226371	a
1.682598	64.400267	a
0.66349	19.044396	a
0.516959	22.175418	a
0.514386	70.405004	a
0.58839	134.752461	a
4.207021	351.258004	a
4.358473	191.316929	a
0.394936	228.437161	a
0.664371	15.601177	a
0.555483	156.678269	a
4.879342	26.477694	a
0.668998	87.548316	a
1513.83165	1298.479651	a
19.516277	36.625863	a
1.130225	33.633605	a
0.232439	125.322952	a
1959.199236	6828.688814	a
10.246557	1770.607378	a
6.77909	468.553289	a
17.549154	1615.598942	a
0.540062	336.845867	a
40.90893	64.502234	a
0.883096	38.904753	a
0.876595	70.547615	a
44.223299	538.126845	a
1.044204	143.373011	a
0.496524	457.694362	a
3142.928333	861.481751	b
180.471772	2670.975838	b
10.24145	1253.117008	b
1.456811	276.530983	b
0.688507	1549.652998	b
3633.046659	1244.080504	b
10029.88386	10000.24815	b
520.742421	1732.087534	b
23.926432	4708.430925	b
302.534675	155.667418	b
0.688346	246.616557	b
15.367324	799.856536	b
6251.247027	1108.91235	b
270.15533	5902.29224	b
142.993533	10000.9523	b
25.62838	542.13942	b
10003.30636	5587.211396	b
23.810349	445.210218	b
10028.97741	10003.7284	b
1.408251	1823.247539	b
87.232802	9632.581557	b
10002.30044	10003.38659	b
2644.616579	10002.28565	b
1.26797	562.603089	b
10003.55686	10003.60447	b
0.941731	10000.9797	b
10000	10000.18984	b
1.929184	1904.455337	b
326.950743	5926.116826	b
1.359575	2238.531849	b
10037.12485	10000.10053	b
0.689911	7607.96861	b
};
\draw [dashed] (0,0) -- (120,120);
	\end{loglogaxis}
\end{tikzpicture}
            \caption{\qIP\ vs \qBB.\vspace{0.3\baselineskip} \quad}    
            \label{fig:genQP_BB}
        \end{subfigure}
        \begin{subfigure}{0.495\textwidth}  
            \centering 
            \begin{tikzpicture}[scale = \myplotscale]
	\begin{loglogaxis}[%
	height = 8.5cm,
	width = 10cm,
	xmin=0,
  xmax=12000,
  ymin=0,
  ymax=12000,
    xlabel={{\tt quadprogIP}},
    ylabel={{\tt \BARON{}}},
	scatter/classes={%
		a={mark=square,black, scale = 0.75},
		b={mark=square,black, scale = 1.2},
		c={mark=square,black, scale = 1.8}}]
	\addplot[scatter,only marks,%
		scatter src=explicit symbolic]%
	table[meta=label] {
x     y      label
6.272094	0.18	a
2.680431	0.02	a
6.223431	0.17	a
3.711449	0.1	a
4.633106	0.08	a
0.198511	0.01	a
0.316398	0.2	a
0.496636	0.22	a
0.191113	0.04	a
0.301215	1.02	a
0.191289	1.13	a
0.199246	0.12	a
0.287471	0.32	a
8.14685	1.46	a
0.431253	0.7	a
0.259599	4.75	a
0.19084	0.02	a
12.180735	10000	b
13.183473	10000.04	b
0.198109	0.02	a
0.198224	0.02	a
0.258417	0.11	a
0.286848	0.1	a
0.197471	0.03	a
0.464387	0.12	a
0.197509	0.02	a
0.264805	0.28	a
0.297753	0.06	a
0.204523	0.01	a
0.195014	0.06	a
0.336474	0.02	a
0.199258	0.03	a
0.400184	0.18	a
0.196603	0.04	a
0.411414	1.02	a
0.278064	1.12	a
0.203168	0.12	a
0.334686	0.32	a
4.104725	0.69	a
5.435408	0.69	a
5.108108	0.72	a
3.915188	0.69	a
8.022685	1.45	a
153.318189	10000	c
0.215762	0.08	a
0.202918	0.13	a
0.250505	0.11	a
0.19602	10000	a
0.256755	0.13	a
0.203145	0.06	a
0.199787	0.1	a
0.200047	0.1	a
0.273755	0.06	a
0.340405	0.14	a
0.270989	0.14	a
0.38285	0.03	a
0.197201	0.01	a
0.283756	0.07	a
0.879105	0.49	a
12.502465	1.25	a
0.2821	0.08	a
0.323705	0.19	a
0.295503	0.12	a
0.192499	0.17	a
0.588008	0.06	a
0.382033	0.06	a
0.295654	0.06	a
0.243686	0.06	a
0.267263	0.12	a
0.268503	0.1	a
0.31459	0.07	a
0.27726	0.08	a
0.256017	0.07	a
0.250592	0.06	a
0.309	0.03	a
0.292584	0.04	a
0.252223	0.04	a
0.191352	0.04	a
0.191123	0.08	a
0.30116	0.02	a
0.337077	0.28	a
0.424575	0.77	a
0.283947	0.31	a
0.558561	0.61	a
4.675576	0.88	a
17.060661	0.62	a
36.991265	0.7	b
125.205971	4.26	b
2.99288	0.82	a
4.536537	0.12	a
0.610013	0.42	a
0.646539	5.36	a
1.540013	0.66	a
8.264153	10000	a
0.389049	1.17	a
1.682598	10000	a
0.66349	253.7	a
0.516959	10000	a
0.514386	1.88	a
0.58839	3.41	a
4.207021	9.49	a
4.358473	1.29	a
0.394936	27.74	a
0.664371	10000	a
0.555483	5.92	a
4.879342	10.16	a
0.668998	0.04	a
1513.83165	2.16	a
19.516277	10000	a
1.130225	0.51	a
0.232439	3.14	a
1959.199236	30.33	a
10.246557	58.28	a
6.77909	10000	a
17.549154	53.5	a
0.540062	10.71	a
40.90893	10000	a
0.883096	25.66	a
0.876595	10000	a
44.223299	10000.01	a
1.044204	151.1	a
0.496524	14.74	a
3142.928333	10000	b
180.471772	10000	b
10.24145	10000	b
1.456811	28.54	b
0.688507	31.2	b
3633.046659	10000.01	b
10029.88386	765.34	b
520.742421	0.09	b
23.926432	218.51	b
302.534675	10000	b
0.688346	191.57	b
15.367324	10000.03	b
6251.247027	10000.01	b
270.15533	314.23	b
142.993533	10000	b
25.62838	10000.01	b
10003.30636	116.32	b
23.810349	10000	b
10028.97741	10000	b
1.408251	10000	b
87.232802	87.06	b
10002.30044	370.45	b
2644.616579	10000.02	b
1.26797	10000	b
10003.55686	1175.22	b
0.941731	35.08	b
10000	10000.02	b
1.929184	5290.3	b
326.950743	10000.01	b
1.359575	3084.81	b
10037.12485	10000.05	b
0.689911	430.75	b
};
\draw [dashed] (-8.5,-10) -- (95,100);
	\end{loglogaxis}
\end{tikzpicture}
            \caption[]%
            {\qIP\ vs \BARON.\vspace{0.3\baselineskip}  \quad}    
            \label{fig:genQP_BA}
        \end{subfigure}
        \begin{subfigure}{0.495\textwidth}   
            \centering 
            \begin{tikzpicture}[scale = \myplotscale]
	\begin{loglogaxis}[%
	height = 8.5cm,
	width = 10cm,
	xmin=0,
  xmax=12000,
  ymin=0,
  ymax=12000,
    xlabel={{\tt quadprogIP}},
    ylabel={{\tt CPLEX}},
	scatter/classes={%
		a={mark=square,black, scale = 0.75},
		b={mark=square,black, scale = 1.2},
		c={mark=square,black, scale = 1.8}}]
	\addplot[scatter,only marks,%
		scatter src=explicit symbolic]%
	table[meta=label] {
x     y      label
6.272094	3.737751	a
2.680431	0.00929	a
6.223431	3.524408	a
3.711449	2.017985	a
4.633106	2.2192	a
0.198511	0.001594	a
0.316398	0.122811	a
0.496636	0.275077	a
0.191113	0.002604	a
0.301215	0.590846	a
0.191289	0.011584	a
0.199246	0.372965	a
0.287471	0.325865	a
8.14685	0.487639	a
0.431253	0.195919	a
0.259599	0.421419	a
0.19084	0.002664	a
12.180735	0.014654	b
13.183473	0.015336	b
0.198109	0.150471	a
0.198224	0.085349	a
0.258417	0.182956	a
0.286848	0.400069	a
0.197471	0.226189	a
0.464387	0.096549	a
0.197509	0.002524	a
0.264805	0.208749	a
0.297753	0.266571	a
0.204523	0.271834	a
0.195014	0.238644	a
0.336474	0.218883	a
0.199258	0.112132	a
0.400184	0.135062	a
0.196603	0.002577	a
0.411414	0.112343	a
0.278064	0.00792	a
0.203168	0.264194	a
0.334686	0.124403	a
4.104725	0.359302	a
5.435408	0.442844	a
5.108108	0.219879	a
3.915188	0.229292	a
8.022685	0.447375	a
153.318189	13445.54623	c
0.215762	0.220919	a
0.202918	0.1516	a
0.250505	0.037848	a
0.19602	0.236114	a
0.256755	0.028441	a
0.203145	0.298705	a
0.199787	0.252112	a
0.200047	0.196568	a
0.273755	0.274709	a
0.340405	0.583759	a
0.270989	0.500927	a
0.38285	0.249514	a
0.197201	0.358374	a
0.283756	0.094534	a
0.879105	0.682915	a
12.502465	0.641891	a
0.2821	0.120671	a
0.323705	0.226139	a
0.295503	0.200346	a
0.192499	0.00479	a
0.588008	0.096156	a
0.382033	0.102075	a
0.295654	0.085491	a
0.243686	0.205828	a
0.267263	0.203342	a
0.268503	0.115148	a
0.31459	0.095063	a
0.27726	0.086783	a
0.256017	0.093457	a
0.250592	0.071478	a
0.309	0.016145	a
0.292584	0.017476	a
0.252223	0.025755	a
0.191352	0.0029	a
0.191123	0.003097	a
0.30116	0.242057	a
0.337077	0.390462	a
0.424575	0.304285	a
0.283947	0.582083	a
0.558561	0.396185	a
4.675576	0.412275	a
17.060661	0.674305	a
36.991265	0.661887	b
125.205971	7.611759	b
2.99288	5.158636	a
4.536537	5.714383	a
0.610013	0.025575	a
0.646539	0.495102	a
1.540013	0.632717	a
8.264153	0.316106	a
0.389049	0.219358	a
1.682598	1.465399	a
0.66349	0.590884	a
0.516959	0.72953	a
0.514386	0.356036	a
0.58839	0.449081	a
4.207021	0.409977	a
4.358473	0.372931	a
0.394936	0.448965	a
0.664371	0.160446	a
0.555483	0.520885	a
4.879342	1.091849	a
0.668998	0.011606	a
1513.83165	0.415729	a
19.516277	1.039526	a
1.130225	0.138693	a
0.232439	0.78251	a
1959.199236	1.592853	a
10.246557	0.98264	a
6.77909	0.619767	a
17.549154	1.897466	a
0.540062	1.114489	a
40.90893	0.370723	a
0.883096	0.655611	a
0.876595	0.27321	a
44.223299	2.15411	a
1.044204	3.068761	a
0.496524	2.842679	a
3142.928333	0.893254	b
180.471772	0.97961	b
10.24145	0.717497	b
1.456811	1.294914	b
0.688507	0.428505	b
3633.046659	1203.98493	b
10029.88386	2.130364	b
520.742421	0.015349	b
23.926432	2.113403	b
302.534675	46.139384	b
0.688346	3.915167	b
15.367324	3.776503	b
6251.247027	512.905863	b
270.15533	9.193145	b
142.993533	37.542279	b
25.62838	26.245546	b
10003.30636	1.177042	b
23.810349	1.882866	b
10028.97741	20.248527	b
1.408251	1.953375	b
87.232802	0.024369	b
10002.30044	2.963123	b
2644.616579	2.73107	b
1.26797	1.97262	b
10003.55686	6.368361	b
0.941731	4.686898	b
10000	135.017745	b
1.929184	24.440686	b
326.950743	94.984709	b
1.359575	34.803928	b
10037.12485	88.563077	b
0.689911	26.634597	b
};
\draw [dashed] (-10,-10) -- (120,120);
	\end{loglogaxis}
\end{tikzpicture}
            \caption[]%
            {\qIP\ vs \cplex.\vspace{0.3\baselineskip}  \quad}    
            \label{fig:genQP_CP}
        \end{subfigure}
        \begin{subfigure}{0.495\textwidth}   
            \centering 
            \includegraphics[width=\textwidth]{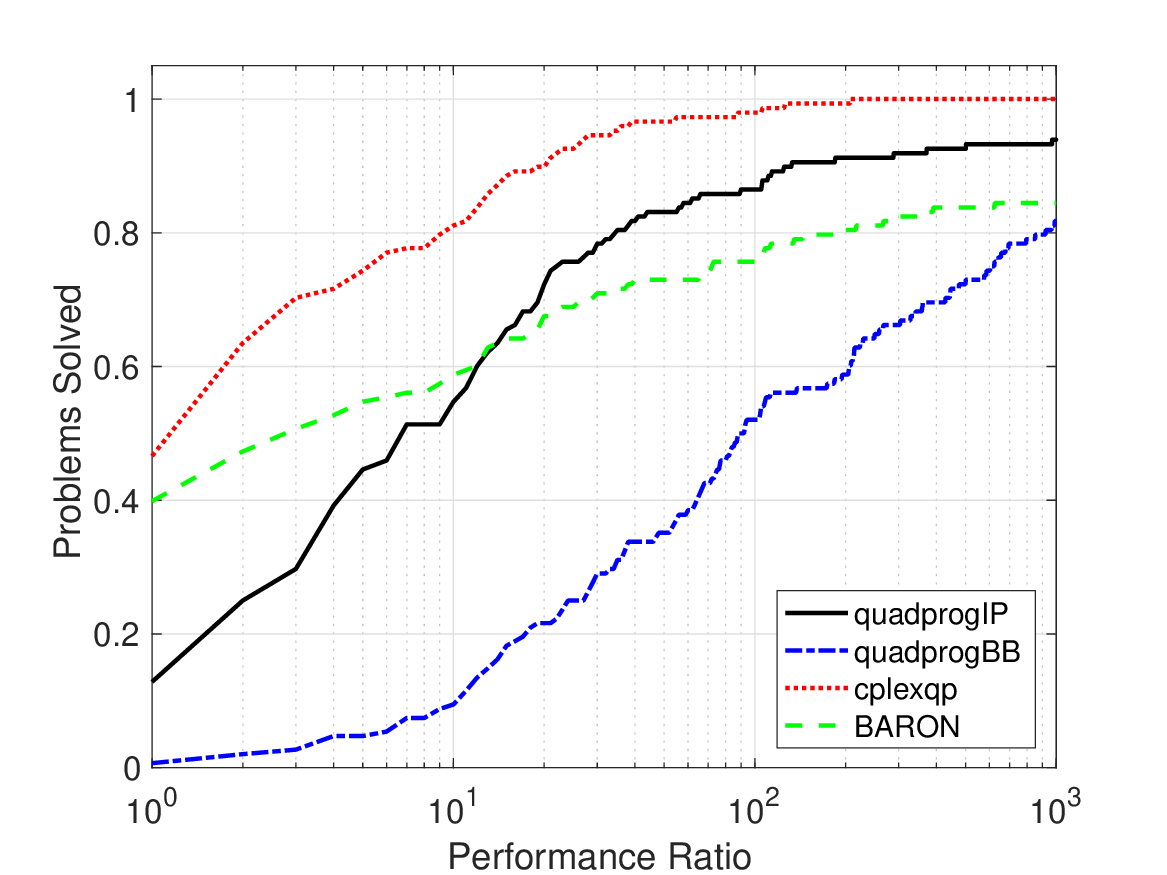}
            \caption[]%
            {Performance profile for CUTEr, Globallib, and RandQP instances.\vspace{0.3\baselineskip}  \quad}    
            \label{fig:genQP_Profile}
        \end{subfigure}
        \caption[]
        {Solution time in seconds of CUTEr, Globallib, and RandQP instances. Size of squares illustrates size of the instance. A square at the maximum value of an axis represents an instances for which the solver in that axis reached maximum running time without solving it.} 
        \label{fig:genQP_ALL}
    \end{figure}
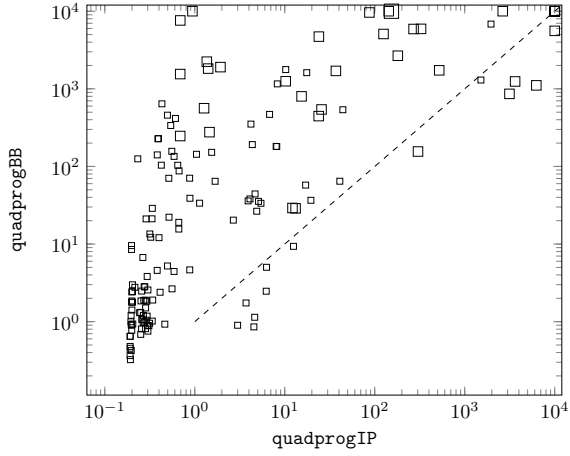
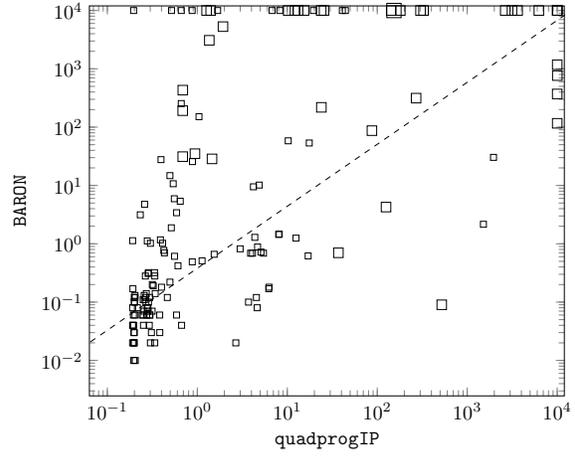
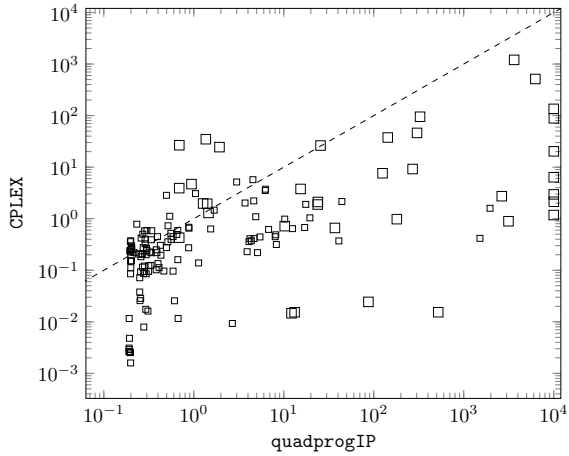
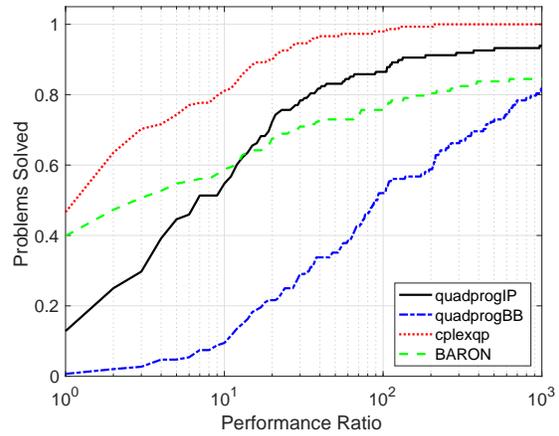

In Figures~\ref{fig:genQP_ALL}, we compare the performance of \qIP\  on the CUTEr, Globallib, and RandQP instances against the other solvers.
As Figure~\ref{fig:genQP_BB} illustrates, except for a few instances, \qIP\  has shorter solution times than \qBB\ on the more general CUTEr, Globallib, and RandQP instances of  QP. Moreover, \qIP\  typically solves
these problems about 10 times faster than \qBB. For these CUTEr, Globallib, and RandQP we find nine (9) instances that are successfully solved by \qIP\  but not by
\qBB\ within the maximum allowed solution time of $10^4$ seconds.

As for \BARON{}, it can be seen from Figures~\ref{fig:genQP_BA} that \qIP is faster on most of the CUTEr, Globallib, and RandQP instances, with \qIP being able to solve a fair number of instances that \BARON is not able to solve within the maximum allowed time of~$10^4$ seconds. On the other hand, \cplex is able to solve most of the CUTEr, Globallib, and RandQP instances faster than \qIP; however, still a number of instances are solved faster than \cplex, and most instances are solved by \qIP in a time no larger than 10 times the solution time of \cplex. Figure~\ref{fig:genQP_Profile} summarizes these results.

\section{Conclusion}
In this paper, we present a new simple and effective approach
for the global solution of  (non-convex) linearly constrained quadratic problems (QP)
by combining the use of the problem's necessary KKT conditions
together with state-of-the-art integer programming solvers. This is done
via a reformulation of the QP as a mixed-integer linear program~(MILP).  We show that
in general, this MILP reformulation can be obtained for QPs with a bounded
primal feasible set via fundamental results related to the solution of
perturbed linear systems of equations \citep[see, e.g.,][]{GuleHR95}. In practice, \qIP is shown to typically outperform by orders of magnitude \qBB, \BARON, and \cplex on standard QPs. For general QPs, \qIP outperforms \qBB, outperforms \BARON in most instances, while \cplex performs the best on these instances. For box-constrained QPs, \qIP has a comparable performance to \qBB and \BARON in small- to medium-scale instances, but is outperformed by these solvers on large-scale instances, while \cplex performs the best on box-constrained QP instances. Also, unlike \qBB, the solution approach proposed here is able to solve QP instances whose dual feasible set is unbounded.
The performance of this methodology on standard QP problems allows for the potential use of this solution
approach as a basis for the solution of {\em copositive programs} \citep[cf.,][]{Dur10}. Which is an interesting direction
of future research work.

The proposed IP formulation of general QPs requires
the computation of certain type of Hoffman bound \citep[see, e.g.,][]{Hoff52} on the system of linear equations defining the problem's
feasible set. Thus, obtaining general and effectively computable bounds of this type is an interesting open question.

We finish by mentioning that a basic implementation of the proposed solution approach referred as \qIP\  is
publicly available at \url{https://github.com/xiawei918/quadprogIP}, together with
pointers to the test instances used in the article for the numerical experiments, and the raw data of all the solution times used to construct the figures throughout the article in the PDF file {\tt raw data.pdf}.

\section*{Acknowledgments}
The authors would like to thank two anonymous referees for the very thoughtful and thorough comments provided on an earlier version of this manuscript.
The work of the first and third author is supported by the NSF CMMI grant 1300193.



\begin{thebibliography}{42}
\expandafter\ifx\csname natexlab\endcsname\relax\def\natexlab#1{#1}\fi
\expandafter\ifx\csname url\endcsname\relax
  \def\url#1{{\tt #1}}\fi
\expandafter\ifx\csname urlprefix\endcsname\relax\def\urlprefix{URL }\fi
\expandafter\ifx\csname urlstyle\endcsname\relax
  \expandafter\ifx\csname doi\endcsname\relax
  \def\doi#1{doi:\discretionary{}{}{}#1}\fi \else
  \expandafter\ifx\csname doi\endcsname\relax
  \def\doi{doi:\discretionary{}{}{}\begingroup \urlstyle{rm}\Url}\fi \fi

\bibitem[{Belotti(2010)}]{Belo10}
Belotti, P. 2010.
\newblock Couenne: a user's manual.
\newblock Tech. rep., Clemson University.
\newblock Available at \url{http://projects.coin-or.
  org/Couenne/browser/trunk/Couenne/doc/couenne-user-manual.pdf}.

\bibitem[{Belotti et~al.(2009)Belotti, Lee, Liberti, Margot, and
  Wachter}]{BeloLLMW09}
Belotti, P., J.~Lee, L.~Liberti, F.~Margot, A.~Wachter. 2009.
\newblock Branching and bounds tightening techiques for non-convex {MINLP}.
\newblock {\it Optimization Methods and Software\/} {\bf 24} 597--634.

\bibitem[{Ben-Tal and Nemirovski(2001)}]{ben-tal01}
Ben-Tal, A., A.~Nemirovski. 2001.
\newblock {\it Lectures on Modern Convex Optimization: Analysis, Algorithms,
  and Engineering Applications\/}.
\newblock MPS-SIAM Series on Optimization, SIAM, Philadelphia, PA.

\bibitem[{Bertsekas(1999)}]{Bert99}
Bertsekas, D. 1999.
\newblock {\it Nonlinear Programming\/}.
\newblock Athena Scientific, Belmont, MA.

\bibitem[{Bomze(1998)}]{Bomze98}
Bomze, I.~M. 1998.
\newblock On standard quadratic optimization problems.
\newblock {\it Journal of Global Optimization\/} {\bf 13} 369--387.

\bibitem[{Bomze et~al.(2010)Bomze, Frommlet, and Locatelli}]{bom2010}
Bomze, I.~M., F.~Frommlet, M.~Locatelli. 2010.
\newblock Copositivity cuts for improving {SDP} bounds on the clique number.
\newblock {\it Mathematical Programming\/} {\bf 124} 13--32.

\bibitem[{Bonami et~al.(2016{\natexlab{a}})Bonami, Lodi, Schweiger, and
  Tramontani}]{bonami2016cut}
Bonami, P., A.~Lodi, J.~Schweiger, A.~Tramontani. 2016{\natexlab{a}}.
\newblock Solving standard quadratic programming by cutting planes.
\newblock Tech. Rep. DS4DM-2016-001, hola.
\newblock Available at
  \url{http://cerc-datascience.polymtl.ca/wp-content/uploads/2016/06/Technical-Report_DS4DM-2016-001-1.pdf}.

\bibitem[{Bonami et~al.(2016{\natexlab{b}})Bonami, G{\"u}nl{\"u}k, and
  Linderoth}]{bonami2016solving}
Bonami, Pierre, Oktay G{\"u}nl{\"u}k, Jeff Linderoth. 2016{\natexlab{b}}.
\newblock Solving box-constrained nonconvex quadratic programs.
\newblock {\it Optimization Online\/} .

\bibitem[{Bundfuss and D{\"u}r(2009)}]{bun2009}
Bundfuss, S., D{\"u}r. 2009.
\newblock An adaptive linear approximation algorithm for copositive programs.
\newblock {\it SIAM Journal on Optimization\/} {\bf 20} 30--53.

\bibitem[{Burer(2009)}]{burer2009copo}
Burer, S. 2009.
\newblock On the copositive representation of binary and continuous nonconvex
  quadratic programs.
\newblock {\it Mathematical Programming\/} {\bf 120} 479--495.

\bibitem[{Burer(2010)}]{BurerDNN}
Burer, S. 2010.
\newblock Optimizing a polyhedral-semidefinite relaxation of completely
  positive programs.
\newblock {\it Mathematical Programming Computation\/} {\bf 2} 1--19.

\bibitem[{Burer and Vandenbussche(2009)}]{burer2009}
Burer, S., D.~Vandenbussche. 2009.
\newblock Globally solving box-constrained non-convex quadratic programs with
  semidefinite-based finite branch-and-bound.
\newblock {\it Computational Optimization and Applications\/} {\bf 43}
  181--195.

\bibitem[{Chen and Burer(2012)}]{chen2012}
Chen, J., S.~Burer. 2012.
\newblock Globally solving nonconvex quadratic programming problems via
  completely positive programming.
\newblock {\it Mathematical Programming Computation\/} {\bf 4} 33--52.

\bibitem[{{CPLEX}(2010)}]{cplex201012}
{CPLEX}, {IBM}~{ILOG}. 2010.
\newblock 12.2 user?s manual.
\newblock {\it ILOG. See ftp://ftp. software. ibm.
  com/software/websphere/ilog/docs/optimization/cplex/ps\_usrmancplex. pdf\/} .

\bibitem[{Dobre and Vera(2015)}]{do2013}
Dobre, C., J.~C. Vera. 2015.
\newblock Exploiting symmetry in copositive programs via semidefinite
  hierarchies.
\newblock {\it Mathematical Programming\/} {\bf 151} 659--680.

\bibitem[{Dong and Anstreicher(2013)}]{dong2013}
Dong, H.~B., K.~Anstreicher. 2013.
\newblock Separating doubly nonnegative and completely positive matrices.
\newblock {\it Mathematical Programming\/} {\bf 137} 131--153.

\bibitem[{D{\"u}r(2010)}]{Dur10}
D{\"u}r, M. 2010.
\newblock Copositive programming -- a survey.
\newblock M.~Diehl, F.~Glineur, E.~Jarlebring, W.~Michiels, eds., {\it Recent
  Advances in Optimization and its Applications in Engineering\/}. Springer,
  3--20.

\bibitem[{Eustaquio et~al.(2008)Eustaquio, Karas, and
  Ribeiro}]{eustaquio2008constraint}
Eustaquio, RODRIGO~G, ELIZABETH~W Karas, ADEMIR~A Ribeiro. 2008.
\newblock Constraint qualifications for nonlinear programming.
\newblock {\it Federal University of Parana\/} .

\bibitem[{Floudas and Visweswaran(1990)}]{flou1990}
Floudas, C.~A., V.~Visweswaran. 1990.
\newblock A global optimization algorithm {(GOP)} for certain classes of
  nonconvex {NLP}s--{I}. {T}heory.
\newblock {\it Computers and Chemical Engineering\/} {\bf 14} 1397--1417.

\bibitem[{Gao(2004)}]{gao2004}
Gao, D.~Y. 2004.
\newblock Canonical duality theory and solutions to constrained nonconvex
  quadratic programming.
\newblock {\it Journal of Global Optimization\/} {\bf 29} 377--399.

\bibitem[{Giannessi and Tomasin(1973)}]{Gian73}
Giannessi, F., E.~Tomasin. 1973.
\newblock Nonconvex quadratic programs, linear complementarity problems, and
  integer linear programs.
\newblock {\it Lecture Notes in Computer Science\/}, vol.~3. Fifth Conference
  on Optimization Techniques, Springer, Berlin Heidelberg New York, 437--449.

\bibitem[{Gould et~al.(2003)Gould, Orban, and Toint}]{gould2003}
Gould, N. I.~M., D.~Orban, P.L. Toint. 2003.
\newblock {CUTE}r and {S}if{D}ec: A constrained and unconstrained testing
  environment, revisited.
\newblock {\it ACM Transactions on Mathematical Software\/} {\bf 29} 373--394.

\bibitem[{G{\"u}ler et~al.(1995)G{\"u}ler, Hoffman, and Rothblum}]{GuleHR95}
G{\"u}ler, O., A.~J. Hoffman, U.~G. Rothblum. 1995.
\newblock Approximations to solutions to systems of linear inequalities.
\newblock {\it SIAM Journal on Matrix Analysis and Applications\/} {\bf 16}
  688--696.

\bibitem[{Hansen et~al.(1993)Hansen, Jaumard, Ruiz, and
  Xiong}]{hansen1993global}
Hansen, Pierre, Brigitte Jaumard, Mich{\`e}le Ruiz, Junjie Xiong. 1993.
\newblock Global minimization of indefinite quadratic functions subject to box
  constraints.
\newblock {\it Naval Research Logistics (NRL)\/} {\bf 40} 373--392.

\bibitem[{Hoffman(2003)}]{Hoff52}
Hoffman, Alan~J. 2003.
\newblock On approximate solutions of systems of linear inequalities.
\newblock {\it Selected Papers Of Alan J Hoffman: With Commentary\/}  174--176.

\bibitem[{Horst et~al.(2000)Horst, Pardalos, and Thoai}]{HorsPT00}
Horst, R., P.~M. Pardalos, N.V. Thoai. 2000.
\newblock {\it Introduction to Global Optimization\/}.
\newblock 2nd ed. Dortrecht: Kluwer.

\bibitem[{Hu et~al.(2012)Hu, Mitchell, Pang, and Yu}]{Pang12}
Hu, J., J.~E. Mitchell, J.~S. Pang, B.~Yu. 2012.
\newblock On linear programs with linear complementarity constraints.
\newblock {\it Journal of Global Optimization\/} {\bf 53} 29--51.

\bibitem[{Kim and Kojima(2001)}]{kim2001}
Kim, S., M.~Kojima. 2001.
\newblock Second order cone programming relaxation of nonconvex quadratic
  optimization problems.
\newblock {\it Optimization Methods and Software\/} {\bf 15} 201--224.

\bibitem[{Kim and Kojima(2003)}]{kim2003}
Kim, S., M.~Kojima. 2003.
\newblock Exact solutions of some nonconvex quadratic optimization problems via
  {SDP} and {SOCP} relaxations.
\newblock {\it Computational Optimization and Applications\/} {\bf 26}
  143--154.

\bibitem[{Mangasarian(1981)}]{Mang81}
Mangasarian, O.~L. 1981.
\newblock A condition number for linear inequalities and linear programs.
\newblock Tech. Rep. 2231, MRC Technical Summary Report.

\bibitem[{Misener and Floudas(2013)}]{misener2013glomiqo}
Misener, Ruth, Christodoulos~A Floudas. 2013.
\newblock Glomiqo: Global mixed-integer quadratic optimizer.
\newblock {\it Journal of Global Optimization\/} {\bf 57} 3--50.

\bibitem[{Motzkin and Straus(1965)}]{mot1965}
Motzkin, T.~S., E.~G. Straus. 1965.
\newblock Maxima for graphs and a new proof of a theorem of {T}ur{\'a}n.
\newblock {\it Canadian Journal of Mathematics\/} {\bf 17} 533--540.

\bibitem[{Nesterov(1998)}]{nesterov1998}
Nesterov, Y. 1998.
\newblock Semidefinite relaxation and nonconvex quadratic optimization.
\newblock {\it Optimization Methods and Software\/} {\bf 9} 141--160.

\bibitem[{Pardalos and Vavasis(1991)}]{PardV91}
Pardalos, P.~M., S.~A. Vavasis. 1991.
\newblock Quadratic programming with one negative eigenvalue is {NP}-hard.
\newblock {\it Journal of Global Optimization\/} {\bf 1} 15--22.

\bibitem[{Pe{\~n}a et~al.(2017)Pe{\~n}a, Vera, and Zuluaga}]{Penahoff17}
Pe{\~n}a, Javier, Juan~C. Vera, Luis~F. Zuluaga. 2017.
\newblock Hoffman bounds and norms of set-valued mappings.
\newblock Tech. rep., Carnegie Mellon University.

\bibitem[{Renegar(2001)}]{Rene01}
Renegar, J. 2001.
\newblock {\it A Mathematical View of Interior-Point Methods in Convex
  Optimization\/}, {\it MPS/SIAM Series on Optimization\/}, vol.~3.
\newblock SIAM.

\bibitem[{Sahinidis(1996)}]{Sahi96}
Sahinidis, N.~V. 1996.
\newblock {BARON}: a general purpose global optimization software package.
\newblock {\it Journal of Global Optimization\/} {\bf 8} 201--205.

\bibitem[{Scozzari and Tardella(2008)}]{scozzari2008clique}
Scozzari, Andrea, Fabio Tardella. 2008.
\newblock A clique algorithm for standard quadratic programming.
\newblock {\it Discrete Applied Mathematics\/} {\bf 156} 2439--2448.

\bibitem[{Sherali and Adams(1994)}]{SherA94}
Sherali, H., W.~Adams. 1994.
\newblock A hierarchy of relaxations for mixed-integer zero-one programming
  problems.
\newblock {\it Discrete Applied Mathematics\/} {\bf 52} 83--106.

\bibitem[{Tawarmalani and Sahinidis(2004)}]{TamaS04}
Tawarmalani, M., N.~V. Sahinidis. 2004.
\newblock Global optimization of mixed-integer nonlinear programs: A
  theoretical and computational study.
\newblock {\it Mathematical Programming\/} {\bf 99} 563--591.

\bibitem[{Vanderbei and Shanno(1999)}]{van1999}
Vanderbei, R.~J., D.~F. Shanno. 1999.
\newblock An interior-point algorithm for nonconvex nonlinear programming.
\newblock {\it Computational Optimization and Applications\/} {\bf 13}
  231--252.

\bibitem[{Zheng and Ng(2004)}]{Ng04}
Zheng, X.~Y., K.~F. Ng. 2004.
\newblock Hoffman's least error bounds for systems of linear inequalities.
\newblock {\it Journal of Global Optimization\/} {\bf 30} 391--403.

\end{thebibliography}
\end{document}